\RequirePackage{ifpdf}
\ifpdf 
\documentclass[pdftex]{sigma}
\else
\documentclass{sigma}
\fi

\numberwithin{equation}{section}
\numberwithin{theorem}{section}
\numberwithin{proposition}{section}
\numberwithin{lemma}{section}
\numberwithin{corollary}{section}
\numberwithin{definition}{section}
\numberwithin{example}{section}
\numberwithin{remark}{section}
\numberwithin{note}{section}

\begin{document}


\renewcommand{\thefootnote}{$\star$}

\renewcommand{\PaperNumber}{094}

\FirstPageHeading

\ShortArticleName{On Tanaka's Prolongation Procedure for Filtered Structures of Constant Type}

\ArticleName{On Tanaka's Prolongation Procedure\\ for Filtered Structures of Constant Type\footnote{This paper is a
contribution to the Special Issue ``\'Elie Cartan and Dif\/ferential Geometry''. The
full collection is available at
\href{http://www.emis.de/journals/SIGMA/Cartan.html}{http://www.emis.de/journals/SIGMA/Cartan.html}}}

\Author{Igor ZELENKO}

\AuthorNameForHeading{I.~Zelenko}

\Address{Department of Mathematics, Texas A$\&$M University,
   College Station, TX 77843-3368, USA}

\Email{\href{mailto:zelenko@math.tamu.edu}{zelenko@math.tamu.edu}}

\ArticleDates{Received June 02, 2009, in f\/inal form September 29, 2009;  Published online October 06, 2009}

\newcommand{\norm}[1]{\left\Vert#1\right\Vert}
\newcommand{\abs}[1]{\left\vert#1\right\vert}
\newcommand{\set}[1]{\left\{#1\right\}}
\newcommand{\Real}{\mathbb R}
\newcommand{\eps}{\varepsilon}
\newcommand{\To}{\longrightarrow}
\newcommand{\BX}{\mathbf{B}(X)}
\newcommand{\A}{\mathcal{A}}
\newcommand{\mg}{\mathfrak g}
\newcommand{\vf}{\varphi}

\Abstract{We present Tanaka's prolongation procedure for f\/iltered structures on manifolds discovered in [Tanaka N.,
{\it J.~Math.\ Kyoto.\ Univ.\/} \textbf{10} (1970), 1--82] in a spirit
  of Singer--Sternberg's  description of  the prolongation of usual $G$-structures
  [Singer I.M., Sternberg S., {\it J.~Analyse Math.\/} {\bf 15} (1965), 1--114;
  Sternberg S.,
Prentice-Hall, Inc., Englewood Clif\/fs, N.J., 1964].
 This approach gives a transparent point of view on the Tanaka constructions avoiding many technicalities of the original Tanaka paper.}


\Classification{58A30; 58A17}

\renewcommand{\thefootnote}{\arabic{footnote}}
\setcounter{footnote}{0}

\section{Introduction}
This note is based on series of lectures given by the author in the Working Geometry Seminar at the Department of Mathematics at Texas A$\&$M University in Spring 2009. The topic is the prolongation procedure for f\/iltered structures on manifolds discovered by Noboru Tanaka in the paper \cite{tan} published in 1970. The Tanaka prolongation procedure is an ingenious ref\/inement  of Cartan's method of equivalence. It provides an ef\/fective algorithm for the construction of canonical frames for f\/iltered structures, and for the calculation of the sharp upper bound of the dimension of their algebras of inf\/initesimal symmetries. This note is by no means a complete survey of the Tanaka theory. For such a survey we refer the reader to \cite{mor}.
 Our goal here is to describe  geometric aspects of Tanaka's prolongation procedure using the language similar to one used by Singer and Sternberg
 in \cite{sinstern} and \cite{stern} for description of the prolongation of the usual $G$-structures.
 We found that it gives a quite natural and transparent point of view on Tanaka's constructions, avoiding many formal def\/initions and technicalities of the original Tanaka paper. We believe this point of view will be useful to anyone who is interested in studying both the main ideas and the details of this fundamental Tanaka construction. We hope that the material of Sections~\ref{section3} and~\ref{section4} will be of interest to experts as well.
 Our language also allows to generalize the Tanaka procedure in several directions, including f\/iltered structures with non-constant and non-fundamental symbols. These generalizations, with applications to the local geometry of distributions, will be given in a separate paper.

\subsection{Statement of the problem}\label{section1.1}

Let $D$ be a rank $l$ distribution on a manifold $M$; that is, a rank $l$ subbundle of the tangent bundle~$TM$.
Two vector distributions
$D_1$ and $D_2$ are called equivalent if there exists a~dif\/feomorphism $F:M\rightarrow M$ such that
$F_*D_1(x)=D_2(F(x))$ for any $x\in M$. Two germs of vector
distributions~$D_1$ and~$D_2$ at the point $x_0\in M$ are
called equivalent, if there exist neighborhoods~$U$ and~$\tilde U$ of~$x_0$ and a dif\/feomorphism $F:U\rightarrow \tilde
U$ such that
\begin{gather*}
F_*D_1=
D_2
, \qquad F(x_0)=x_0.
\end{gather*}

 The general question is: When are two germs of distributions
 equivalent?

\subsection{Weak derived f\/lags and symbols of distributions}\label{section1.2}

Taking Lie brackets of vector f\/ields tangent to a distribution $D$ (i.e.\ sections of $D$)   one can def\/ine a f\/iltration $D^{-1}\subset D^{-2}\subset\cdots$
of the tangent bundle, called a \emph{weak derived flag}  or a \emph{small flag $($of $D)$}. More precisely, set $D=D^{-1}$ and def\/ine recursively
$D^{-j}=D^{-j+1}+[D,D^{-j+1}]$, $j>1$.
Let  $X_1,\ldots X_l$ be $l$ vector f\/ields constituting a local basis of a distribution~$D$, i.e.\ $D= {\rm span}\{X_1, \ldots, X_l\}$
in some open set in $M$. Then $D^{-j}(x)$ is the linear span of all iterated Lie brackets of these vector f\/ields, of length not greater than  $j$,  evaluated at a point~$x$.
A~distribution~$D$ is called \emph{bracket-generating} (or \emph{completely nonholonomic}) if for any $x$ there exists $\mu(x)\in\mathbb N$ such that $D^{-\mu(x)}(x)=T_x M$. The number $\mu(x)$ is called the \emph{degree of nonholonomy} of $D$ at a point~$x$.
A distribution $D$ is called \emph{regular} if for all $j<0$, the dimensions of subspaces $D^j(x)$ are independent of the point $x$.
From now on we assume that $D$ is regular bracket-generating distribution with degree of nonholonomy $\mu$.
Let $\mg^{-1}(x)\stackrel{\text{def}}{=}D^{-1}(x)$ and $\mg^{j}(x)\stackrel{\text{def}}{=}D^{j}(x)/D^{j+1}(x)$ for $j<-1$. Consider the graded space
\begin{gather*}
\mathfrak{m}(x)=\bigoplus_{j=-\mu}^{-1}\mg^j(x),
\end{gather*}
corresponding to the f\/iltration
\begin{gather*}
D(x)=D^{-1}(x)\subset D^{-2}(x)\subset\cdots\subset D^{-\mu+1}(x)
\subset D^{-\mu}(x)=T_xM.
\end{gather*}
 This space is endowed naturally with the structure of a graded nilpotent Lie algebra, generated by
$\mg^{-1}(x)$. Indeed, let $\mathfrak p_j:D^j(x)\mapsto \mg^j(x)$ be the canonical projection to a factor space. Take $Y_1\in\mg^i(x)$ and $Y_2\in \mg^j(x)$. To def\/ine the Lie bracket $[Y_1,Y_2]$ take a local section $\widetilde Y_1$ of the distribution $D^i$ and
a local section $\widetilde Y_2$ of the distribution $D^j$ such that $\mathfrak p_i\bigl(\widetilde Y_1(x)\bigr) =Y_1$
and $\mathfrak p_j\bigl(\widetilde Y_2(x)\bigr)=Y_2$. It is clear that $[Y_1,Y_2]\in\mg^{i+j}(x)$. Put
\begin{gather}
\label{Liebrackets}
[Y_1,Y_2]\stackrel{\text{def}}{=}\mathfrak p_{i+j}\bigl([\widetilde Y_1,\widetilde Y_2](x)\bigr).
\end{gather}
It is easy to see that the right-hand side  of \eqref{Liebrackets} does not depend on the choice of sections $\widetilde Y_1$ and
$\widetilde Y_2$. Besides, $\mg^{-1}(x)$ generates the whole algebra $\mathfrak{m}(x)$.
A graded Lie algebra satisfying the last property is called \emph{fundamental}.
The graded nilpotent Lie algebra $\mathfrak{m}(x)$ is called the \emph {symbol of the distribution $D$ at the point $x$}.

Fix a fundamental graded nilpotent Lie algebra $\mathfrak{m}=\displaystyle{\bigoplus_{i=-\mu}^{-1} \mg^i}$.
A distribution $D$ is said to be of \emph{constant symbol $\mathfrak{m}$} or of \emph{constant type $\mathfrak{m}$} if for any
 $x$ the symbol $\mathfrak{m}(x)$ is isomorphic to $\mathfrak{m}$ as a nilpotent graded Lie algebra.
In general this assumption is quite restrictive.
For example, in the case of rank two distributions on manifolds with $\dim\,M\geq 9$,
symbol algebras  depend on continuous parameters, which implies that generic rank 2 distributions in these dimensions do not have a
constant symbol.
For rank 3 distributions with $\dim D^{-2}=6$  the same holds in the case $\dim M=7$ as was shown in \cite{kuz}.
Following Tanaka, and for simplicity of presentation, we consider here distributions of constant type $\mathfrak{m}$ only.
One can construct the \emph{flat distribution~$D_{\mathfrak{m}}$ of constant type $\mathfrak{m}$}.
For this let $M(\mathfrak{m})$ be the simply connected Lie group with the
Lie algebra~$\mathfrak{m}$ and let $e$ be its identity. Then $D_\mathfrak{m}$ is the left invariant distribution on $M(\mathfrak{m})$ such that $D_{\mathfrak{m}}(e)=\mg^{-1}$.

\subsection[The bundle $P^0(\mathfrak{m})$ and its reductions]{The bundle $\boldsymbol{P^0(\mathfrak{m})}$ and its reductions}\label{section1.3}

\looseness=-1
To a distribution of type $\mathfrak{m}$ one can assign a principal bundle in the following way. Let $G^0(\mathfrak{m})$ be the group of automorphisms
of the graded Lie algebra $\mathfrak{m}$; that is, the group of all automorphisms~$A$ of the linear space $\mathfrak{m}$ preserving both the Lie brackets ($A([v,w])=[A(v),A(w)]$ for any $v,w\in \mathfrak{m}$) and the grading ($A (\mg^i)=\mg^i$ for any $i<0$).
Let
$P^0(\mathfrak m)$ be the set of all pairs $(x,\vf)$, where  $x\in M$ and  $\vf:\mathfrak{m}\to\mathfrak m(x)$
is an isomorphism of the graded Lie algebras~$\mathfrak {m}$ and~$\mathfrak m(x)$.
 Then $P^0(\mathfrak m)$ is a principal $G^0(\mathfrak m)$-bundle over $M$. The right action $R_A$ of an automorphism $A\in G^0(\mathfrak{m})$ is as follows:  $R_A$ sends $(x,\vf)\in P^0(\mathfrak m)$ to $(x,\vf\circ A)$, or shortly $(x,\vf)\cdot R_A=(x,\vf\circ A)$. Note that since $\mg^{-1}$ generates
 $\mathfrak m$, the group $G^0(\mathfrak{m})$ can be identif\/ied with a subgroup of $\text{GL}(\mg^{-1})$. By the same reason a point $(x,\vf)\in P^0(\mathfrak m)$ of a f\/iber of $P^0(\mathfrak m)$ is uniquely def\/ined by $\vf|_{\mg^{-1}}$. So one can identify
$P^0(\mathfrak m)$ with the set of pairs $(x,\psi)$, where  $x\in M $ and $\psi:\mg^{-1}\to D(x)$  can be extended to an automorphism of the graded Lie algebras $\mathfrak {m}$ and $\mathfrak m(x)$.
Speaking informally, $P^0(\mathfrak m)$ can be seen as a $G^0(\mathfrak{m})-$reduction of the bundle of all frames of the distribution $D$.
Besides, the Lie algebra $\mg^0(\mathfrak m)$ is the algebra of all derivations $a$ of $\mathfrak m$, preserving the grading
(i.e.\ $a \mg^i\subset \mg^i$ for all $i<0$).

Additional structures on distributions can be encoded by reductions of the bundle $P^0(\mathfrak m)$.
More precisely, let $G^0$  be a Lie subgroup of $G^0(\mathfrak{m})$ and let $P^0$ be a principal $G^0$-bundle, which is a reduction of the bundle $P^0(\mathfrak{m})$.
Since $\mg^0$ is a subalgebra of the algebra of derivations of $\mathfrak m$ preserving the grading, the subspace
$\mathfrak m\oplus \mg^0$ is endowed with the natural structure of a graded Lie algebra. For this we only need to def\/ine
brackets $[f,v]$
for $f\in \mg^0$ and $v\in \mathfrak m$, because $\mathfrak m$ and $\mg^0$ are already Lie algebras.
Set
$[f,v]\stackrel{\text{def}}{=} f(v)$.
The bundle $P^0$ is called a \emph{structure of constant type $(\mathfrak m, \mg^0)$}.
Let, as before, $D_{\mathfrak m}$ be the left invariant distribution on $M(\mathfrak{m})$ such that $D_{\mathfrak m}(e)=\mg^{-1}$. Denote by $L_x$ the left translation on $M(\mathfrak m)$ by an element $x$. Finally, let $P^0(\mathfrak m, \mathfrak g^0)$
be the set of all pairs $(x,\vf)$,
where  $x\in M(\mathfrak m)$ and  $\vf:\mathfrak{m}\to\mathfrak m(x)$ is an isomorphism of the graded Lie algebras $\mathfrak {m}$ and $\mathfrak m(x)$ such that
$(L_{x^{-1}})_*\vf\in G^0$. The bundle $P^0(\mathfrak m, \mathfrak g^0)$ is called \emph{the flat structure of constant type $(\mathfrak m, \mg^0)$}.
Let us give some examples.

{\bf Example 1. ${\boldsymbol G}$-structures.}
Assume that $D=TM$. So $\mathfrak m=\mg^{-1}$
is abelian, $G^0(\mathfrak m)={\rm GL}(\mathfrak m)$, and $P^0(\mathfrak m)$ coincides with the bundle $\mathcal F(M)$ of all frames on $M$. In this case $P^0$ is nothing but a usual $G^0$-structure.

{\bf Example 2. Contact distributions.} Let $D$ be the contact distribution in $\mathbb R^{2n+1}$. Its symbol
$\mathfrak m_{\text{cont},n}$ is isomorphic to the Heisenberg algebra $\eta_{2n+1}$ with grading $\mg^{-1}\oplus\mg^{-2}$, where $\mg^{-2}$ is the center of
$\eta_{2n+1}$. Obviously, a skew-symmetric form $\Omega$ is well def\/ined on $\mg^{-1}$, up to a multiplication by a nonzero constant.
The group $G^0(\mathfrak m_{\text{cont},n})$ of automorphisms of $\mathfrak m_{\text{cont},n}$ is isomorphic to the group $\text{CSP}(\mg^{-1})$ of conformal symplectic transformations of $\mg^{-1}$, i.e. transformations preserving the form $\Omega$, up to a multiplication by a nonzero constant.

{\bf Example 3. Maximally nonholonomic rank 2 distributions in $\boldsymbol{\mathbb R^5}$.} Let $D$ be a rank~2 distribution in $\mathbb R^5$
with degree of nonholonomy equal to $3$ at every point. Such distributions were treated by \'E.~Cartan in his famous work \cite{cartan}. In this case $\dim D^{-2}\equiv 3$ and $\dim D^{-3}\equiv 5$. The symbol at any point is isomorphic to the Lie algebra $\mathfrak m_{(2,5)}$ generated by $X_1$, $X_2$, $X_3$, $X_4,$ and  $X_5$ with the following nonzero products: $[X_1,X_2]=X_3$, $[X_1, X_3]=X_4$, and $[X_2, X_3]=X_5$. The grading is given as follows:
\[
\mg^{-1}=\langle X_1, X_2\rangle,\qquad \mg^{-2}=\langle X_3 \rangle,\qquad \mg^{-3}=\langle X_4, X_5\rangle,
\]
where $\langle Y_1,\ldots, Y_k\rangle$ denotes the linear span of vectors $Y_1,\ldots,Y_k$.
Since $\mathfrak m_{(2,5)}$ is a free nilpotent Lie algebra with two generators $X_1$ and $X_2$, its group of automorphism is equivalent to $\text{GL}(\mg^{-1})$.

\looseness=-1
{\bf Example 4. Sub-Riemannian structures of constant type} (see also \cite{morsR}).
Assume that each space $D(x)$ is endowed with an Euclidean structure $Q_x$ depending smoothly on $x$.
In this situation
the pair $(D,Q)$ def\/ines a sub-Riemannian structure on a manifold $M$.
Recall that $\mg^{-1}(x)=D(x)$. This motivates  the following def\/inition:
A pair $ \bigl(\mathfrak m, \mathfrak Q)$, where $\mathfrak m=\displaystyle{\bigoplus_{j=-\mu}^{-1} \mg^{j}}$ is a fundamental graded Lie algebra and $\mathfrak Q$ is an Euclidean structure on $\mg^{-1}$, is called a \emph{sub-Riemannian symbol}. Two sub-Riemannian symbols $(\mathfrak m, \mathfrak Q)$ and    $(\tilde{\mathfrak m}, \widetilde {\mathfrak Q})$ are isomorphic if there exists a map $\vf: \mathfrak m\to \tilde{\mathfrak m}$, which is an isomorphism of the graded Lie algebras $\mathfrak m$ and $\tilde{\mathfrak m}$, preserving the Euclidean structures $\mathfrak Q$ and $\tilde{\mathfrak Q}$ (i.e. such that $\widetilde Q \bigl(\vf (v_1), \vf (v_2)\bigr)=Q(v_1, v_2)$ for any $v_1$ and $v_2$ in $\mg^{-1}$). Fix a sub-Riemannian symbol $(\mathfrak m,\mathfrak Q)$. A sub-Riemannian structure $(D,Q)$  is said to be of \emph{constant type $(\mathfrak{m},\mathfrak Q)$}, if for every $x$ the sub-Riemannian symbol $(\mathfrak{m}(x), Q_x)$ is isomorphic to $(\mathfrak{m}, \mathfrak Q)$.

It may happen that a sub-Riemannian structure does not have a constant symbol even if the distribution does.
 Such a situation occurs already in the case of the contact distribution on $\mathbb R^{2n+1}$ for $n>1$ (see Example 2 above).
 As was mentioned above,  in this case a  skew-symmetric form~$\Omega$ is well def\/ined on $\mg^{-1}$, up to a multiplication by a nonzero constant.
 If in addition a Euclidean structure $Q$ is given on $\mg^{-1}$, then a skew-symmetric endomorphism
 $J$ of~$\mg^{-1}$ is well def\/ined, up to a multiplication by a nonzero constant, by $\Omega(v_1,v_2)=Q(J v_1,v_2)$.
 Take $0<\beta_1\leq \cdots\leq\beta_n$ so that $\{\pm \beta_1 i, \ldots,\pm\beta_n i\}$ is the set of
 the eigenvalues of
 $J$. Then a sub-Riemannian symbol with $\mathfrak m= \mathfrak m_{\text{cont},n}$ is determined uniquely (up to an isomorphism)  by a point $[\beta_1:\beta_2:\ldots:\beta_n]$ of the projective space $\mathbb {R P}^{n-1}$.

Let $(D,Q)$  be a sub-Riemannian structure of constant type $(\mathfrak{m},\mathfrak Q)$
and $G^0(\mathfrak{m},\mathfrak Q)\subset G^0(\mathfrak m)$ be the group of automorphisms of a sub-Riemannian symbol $(\mathfrak{m}, \mathfrak Q)$.
Let
$P^0(\mathfrak{m},\mathfrak Q)$ be the set of all pairs $(x,\vf)$, where  $x\in M$ and  $\vf:\mathfrak{m}\to\mathfrak m(x)$ is an isomorhism of sub-Riemannian symbols $\bigl(\mathfrak {m},\mathfrak Q\bigr)$ and $\bigl(\mathfrak m(x), Q_x\bigr)$. Obviously, the bundle $P^0(\mathfrak{m},\mathfrak Q)$ is a reduction of $P^0(\mathfrak{m})$ with the structure group $G^0(\mathfrak{m},\mathfrak Q)$.

\looseness=1
{\bf Example 5. Second order ordinary dif\/ferential equations up to point transformations.} Assume that  $D$ is a contact distribution on a $3$-dimensional manifold
endowed with two distinguished transversal line sub-distributions  $L_1$ and $L_2$. Such structures appear in the study of second order ordinary dif\/ferential equations $y''=F(t,y,y')$ modulo point transformations. Indeed, let $J^i(\mathbb R,\mathbb R)$ be the space of $i$-jets of mappings from $\mathbb R$ to $\mathbb R$. As the distribution~$D$ we take the standard contact distribution on $J^1(\mathbb R,\mathbb R)$. In the standard coordinates $(t,y,p)$ on~$J^1(\mathbb R,\mathbb R)$ this distribution is given by the Pfaf\/f\/ian equation $dy-pdt=0$. The natural lifts to~$J^1$ of solutions of the dif\/ferential equation  form the
$1$-foliation tangent to~$D$.  The tangent lines to this foliation def\/ine  the sub-distribution~$L_1$. In the
coordinates $(t,y,p)$ the sub-distribution~$L_1$ is generated by the vector f\/ield $\frac{\partial}{\partial t}+p\frac{\partial}{\partial y}+F(t,y,p)\frac{\partial}{\partial p}$. Finally, consider the natural bundle
$J^1(\mathbb R,\mathbb R)\rightarrow J^0(\mathbb R,\mathbb R)$ and let $L_2$ be the distribution of the tangent lines to the f\/ibers. The sub-distribution $L_2$ is generated by the vector f\/ield $\frac{\partial}{\partial p}$. The triple $(D,L_1, L_2)$ is called the \emph{pseudo-product structure associated with the second order ordinary differential equation}.
Two second order dif\/ferential equations are equivalent with respect to the group of point transformations if and only if there is a dif\/feomorphism of $J^{1}(\mathbb R,\mathbb R)$ sending  the pseudo-product structure associated with one of them to the pseudo-product structure associated with the other one. This equivalence problem was treated by \'E.~Cartan in~\cite{cartproj} and earlier by A.~Tresse in~\cite{tresse1} and~\cite{tresse2}. The symbol of the distribution is $\mathfrak m_{\text{cont},1}\sim \eta_3$ (see Example~2 above) and the plane $\mg^{-1}$ is endowed with two distinguished transversal lines. This additional structure is encoded by the subgroup~$G^0$ of the group $G^0(\mathfrak m_{\text{cont},1})$ preserving each of these lines.

Another important class of geometric structures that can be encoded in this way are  $CR$-structures (see \S~10 of~\cite{tan}
for more details).

\subsection{Algebraic and geometric Tanaka prolongations}\label{section1.4}

In \cite{tan} Tanaka solves the equivalence problem for structures of constant type $(\mathfrak m, \mg^0)$.
Two of Tanaka's main constructions are the algebraic prolongation of the algebra $\mathfrak m+\mg^0$, and the geometric prolongation of structures of type $(\mathfrak m, \mg^0)$, imitated by the algebraic prolongation.

First he def\/ines a graded Lie algebra,
which is in essence the maximal (nondegenerated) graded Lie algebra, containing the graded Lie algebra
$\displaystyle{\bigoplus_{i\leq 0}\mg^i}$ as its non-positive part. More precisely, Tanaka constructs
a graded Lie algebra $\mg(\mathfrak m, \mg^0)= \displaystyle{\bigoplus_{i\in\mathbb Z}\mg^i(\mathfrak m,\mg^0)}$, satisfying the following three conditions:
\begin{enumerate}\itemsep=0pt
\item $\mg^i(\mathfrak m,\mg^0)=\mg^i$ for all $i\leq 0$;
\item if $X\in \mg^i(\mathfrak m,\mg^0)$ with $i>0$ satisf\/ies $[X, \mg^{-1}]=0$, then $X=0$;
\item $\mg(\mathfrak m, \mg^0)$ is the maximal graded Lie algebra, satisfying Properties 1 and 2.
\end{enumerate}
This graded Lie algebra $\mg(\mathfrak m, \mg^0)$ is called  the \emph{algebraic universal prolongation}  of the graded Lie algebra $\mathfrak m\oplus \mg^0$.
An explicit realization of the algebra $\mg(\mathfrak m, \mg^0)$  will be described later in Section~\ref{section4}.
It turns out (\cite[\S~6]{tan}, \cite[\S~2]{yam}) that the Lie algebra of inf\/initesimal symmetries of the f\/lat structure of type $(\mathfrak m, \mg^0)$ can be described in terms of $\mg(\mathfrak m, \mg^0)$. If $\dim \mg(\mathfrak m, \mg^0)$ is f\/inite (which is equivalent to the existence of $l>0$ such that $\mg^l(\mathfrak m,\mg^0)=0$), then the algebra of inf\/initesimal symmetries is  isomorphic to $\mg(\mathfrak m, \mg^0)$. The analogous formulation in the case  when
$\mg(\mathfrak m, \mg^0)$ is inf\/inite dimensional may be found in \cite[\S~6]{tan}.

 Furthermore for a structure $P^0$ of type $(\mathfrak m, \mg^0)$,  Tanaka constructs a sequence of bundles $\{P^i\}_{i\in\mathbb N}$, where $P^i$ is a principal bundle over $P^{i-1}$ with an abelian structure group  of dimension equal to $\dim \mg^i(\mathfrak m,\mg^0)$. In general $P^i$ is not a frame bundle. This is the case only for $\mathfrak m=\mg^{-1}$; that is, for $G$-structures. But if $\dim \mg(\mathfrak m, \mg^0)$ is f\/inite or, equivalently, if there exists $l\geq 0$ such that $\mg^{l+1}(\mathfrak m,\mg^0)=0$,
then the bundle $P^{l+\mu}$ is an $e$-structure over $P^{l+\mu-1}$, i.e. $P^{l+\mu-1}$ is endowed with a canonical frame (a structure of absolute parallelism).
Note that all $P^i$ with $i\geq l$ are identif\/ied one with each other by the canonical projections (which are dif\/feomorphisms in that case).
Hence,
\emph{$P^{l}$  is endowed with a canonical frame}. Once a canonical frame is constructed the equivalence problem for structures of type $(\mathfrak m, \mg^0)$ is in essence solved. Moreover,  $\dim \mg(\mathfrak m, \mg^0)$ gives the sharp upper bound for the dimension of the algebra of inf\/initesimal symmetries of such
structures.

By Tanaka's geometric prolongation we mean his construction of the sequence of bundles $\{P^i\}_{i\in\mathbb N}$.
In this note we mainly concentrate on a description of this geometric prolongation using a language dif\/ferent from Tanaka's original one. In Section~\ref{section2} we review the
prolongation of usual $G$-structures in the language of Singer and Sternberg.
We do this in order to prepare the reader for the next section, where the f\/irst Tanaka geometric prolongation
is given in a~completely analogous way. We believe that after reading Section~\ref{section3}
the reader will already have an idea how to proceed with the higher order
Tanaka prolongations so that technicalities of Section~\ref{section4} can be easily overcome.

\section[Review of prolongation of $G$-structures]{Review of prolongation of $\boldsymbol{G}$-structures}\label{section2}

Before treating the general case we review the prolongation procedure for structures with \mbox{$\mathfrak m=\mg^{-1}$}, i.e.\ for
usual $G$-structures.
We follow \cite{sinstern} and \cite{stern}.
Let $\Pi_0:P^0\to M$ be the canonical projection and $V(\lambda)\subset T_\lambda P^0$ the tangent space at $\lambda$
to the f\/iber of $P^0$ over the point $\Pi_0(\lambda)$. The subspace $V(\lambda)$ is also called the \emph{ vertical
subspace of $T_\lambda P^0$}. Actually,
\begin{gather}
\label{vert}
V(\lambda)=\ker (\Pi_0)_*(\lambda).
\end{gather}
 Recall that the space $V(\lambda)$ can be identif\/ied with the Lie algebra $\mg^0$ of $G^0$.
 The identif\/ication $I_\lambda:\mg^0\rightarrow V(\lambda)$ sends $X\in \mg^0$ to $\frac{d}{dt}\bigl(\lambda\cdot R_{e^{tX}} \bigr)|_{t=0}$,
 where $e^{tX}$ is the one-parametric subgroup generated by $X$. Recall also that an Ehresmann connection on the bundle $P^0$ is
 a distribution~$H$ on~$P^0$ such that
\begin{gather}
\label{Ehres}
T_\lambda P^0=V(\lambda)\oplus H(\lambda)\qquad \forall\, \lambda\in P^0.
\end{gather}
A subspace $H(\lambda)$, satisfying \eqref{Ehres}, is a \emph{horizontal subspace of $T_\lambda P^0$}.

Once an Ehresmann connection $H$ and a basis in the space $\mg^{-1}\oplus\mg^0$ are f\/ixed, the bund\-le~$P^0$ is endowed with a frame
in a canonical way. Indeed, let $\lambda=(x,\vf)\in P^0$. Then $\vf\in {\rm Hom} (\mg^{-1}, T_xM)$.
By \eqref{vert} and \eqref{Ehres} the restriction $(\Pi_0)_*|_{H(\lambda)}$ of the map $(\Pi_0)_*$ to the subspace~$H(\lambda)$
is an isomorphism between $H(\lambda)$ and $T_{\Pi_0(\lambda)}M$. Def\/ine the map
$\vf^{H(\lambda)}:\mg^{-1}\oplus\mg^0\rightarrow T_\lambda P^0$ as follows:
\begin{gather}
 \vf^{H(\lambda)}|_{\mg^{-1}}=\bigl((\Pi_0)_*|_{H(\lambda)}\bigr)^{-1}\circ\vf,\nonumber\\
 \vf^{H(\lambda)}|_{\mg^0}=I_\lambda.\label{vfH0}
\end{gather}
If one f\/ixes a basis in $\mg^{-1}\oplus\mg^0$, then the images of this basis under the maps $\vf^{H(\lambda)}$ def\/ine the frame (the structure of the absolute parallelism) on~$P^0$.

The question is whether an Ehresmann connection can be chosen canonically. To answer this question, f\/irst one introduces
a special $\mg^{-1}$-valued $1$-form $\omega$ on $P^0$ as follows: $\omega(Y)=\vf^{-1}\circ (\Pi_0)_*(Y)$ for any $\lambda=(x,\vf)\in P^0$ and $Y\in T_\lambda P^0$. This $1$-form is called the \emph{soldering $($tautological, fundamental$)$ form} of the $G^0$-structure
$P^0$. Further, f\/ixing  again a point $\lambda=(x,\vf)\in P^0$, one def\/ines a \emph{structure function $($a torsion$)$} $C_H\in \text {Hom}(\mg^{-1}\wedge\mg^{-1},\mg^{-1})$ of a horizontal subspace $H$ of $T_\lambda P^0$, as follows:
\begin{gather*}
\forall\, v_1,v_2\in \mg^{-1}\qquad C_H(v_1,v_2)=-d\omega\bigl(\vf^H(v_1),\vf^H(v_2)\bigr),
\end{gather*}
where $\vf^H$ is def\/ined by \eqref{vfH0}. Equivalently,
\begin{gather*}
C_H(v_1,v_2)=\omega\bigl([Y_1,Y_2](\lambda)\bigl)
\end{gather*}
for any vector f\/ields $Y_1$ and $Y_2$ such that $\omega(Y_i)\equiv v_i$ and $\vf^H(v_i)=Y_i(\lambda)$, $i=1,2$.
Speaking informally, the structure function $C_H$ encodes all information
about horizontal parts at $\lambda$ of Lie brackets of vector f\/ields which are horizontal at $\lambda$ w.r.t.\
the splitting \eqref{Ehres} (with $H(\lambda)$ replaced by $H$).

We now take another horizontal subspace $\widetilde H$ of $T_\lambda P^0$ and compare the structure
functions~$C_{H}$ and $C_{\widetilde{H}}$.
By construction, for any vector $v\in\mg ^{-1}$ the vector $\vf^{\widetilde{H}}(v)-\vf^{H}(v)$ belongs to
$V(\lambda)$
($\sim \mg^0$).
Let
\begin{gather*}
f_{H\widetilde H}(v)\stackrel{\text{def}}{=}
I_\lambda^{-1}\big(\vf^{\widetilde{H}}(v)-\vf^{H}(v)\big).
\end{gather*}

Then $f_{H\widetilde H}\in {\rm Hom}(\mg^{-1},\mg^0)$. In the opposite direction, it is
clear that for any $f\in {\rm Hom}(\mg^{-1},\mg^{0})$ there exists a horizontal subspace $\widetilde H$
such that  $f=f_{H\widetilde H}$.
The map
\begin{gather*}
\partial: \ {\rm Hom}(\mg^{-1},\mg^{0})\rightarrow
 {\rm Hom}\big(\mg^{-1}\wedge\mg^{-1},\mg^{-1}\big),
\end{gather*}
 def\/ined  by
\begin{gather}
\label{Spencer}
\partial f(v_1,v_2)=f(v_1)v_2- f(v_2)v_1
=[f(v_1),v_2]+[v_1,f(v_2)]
\end{gather}
is called the \emph{Spencer operator}\footnote{In \cite{stern} this operator is called the antisymmetrization operator, but we prefer to call it the Spencer operator, because, after certain intepretation of the spaces ${\rm Hom}(\mg^{-1},\mg^{0})$ and  ${\rm Hom}(\mg^{-1}\wedge\mg^{-1},\mg^{-1})$, this operator can be identif\/ied with an appropriate $\delta$-operator introduced by Spencer in \cite{spenc} for the study of overdetermined systems of partial dif\/ferential equations.  Indeed, since $\mg^0$ is a~subspace of $\mathfrak{gl}(\mg^{-1})$, the space  ${\rm Hom}(\mg^{-1},\mg^{0})$ can be seen as a~subspace of the space of $\mg^{-1}$-valued one-forms on $\mg^{-1}$ with linear coef\/f\/icients, while ${\rm Hom}(\mg^{-1}\wedge\mg^{-1},\mg^{-1})$ can be seen as the space of $\mg^{-1}$-valued two-forms on $\mg^{-1}$
with constant coef\/f\/icients. Then the operator $\partial$ def\/ined  by \eqref{Spencer} coincides with the restriction to  ${\rm Hom}(\mg^{-1},\mg^{0})$ of the exterior dif\/ferential acting between the above-mentioned spaces of one-forms and two-forms, i.e.\ with the corresponding  Spencer $\delta$-operator.}.
By direct computations (\cite[p.~42]{sinstern}, \cite[p.~317]{stern}, or the proof of more general statement in
Proposition \ref{prop1} below)
one obtains the following identity
\begin{gather*}
C_{\widetilde {H}}=C_{H}+\partial f_{H \widetilde{H}}.
\end{gather*}


Now f\/ix a subspace
\begin{gather*}
\mathcal N\subset
 {\rm Hom}\big(\mg^{-1}\wedge\mg^{-1},\mg^{-1}\big)
 \end{gather*}
 complementary to  $\text{Im}\, \partial$, so that
\begin{gather*}
 {\rm Hom}\big(\mg^{-1}\wedge\mg^{-1},\mg^{-1}\big)= \text{Im} \,\partial\oplus\mathcal N.
 \end{gather*}
 Speaking informally, the subspace $\mathcal N$ def\/ines the normalization conditions for
 the f\/irst prolongation.
 The \emph{first prolongation of $P^0$} is the following bundle $(P^0)^{(1)}$ over $P^0$:
 \begin{gather*}
\big(P^0\big)^{(1)} =\big\{(\lambda, H):\lambda\in P^0, H \text{ is a horizontal subspace of } T_\lambda P^0 \text{ with }C_H\in \mathcal N \big\}.
 \end{gather*}
Alternatively,
\begin{gather*}
 \big(P^0\big)^{(1)}=\big\{(\lambda,\vf^H):\lambda\in P^0, H \text{ is a horizontal subspace of } T_\lambda P^0 \text { with } C_H\in \mathcal N\big\}.
 \end{gather*}
 In other words, the f\/iber of $(P^0)^{(1)}$ over a point $\lambda\in P^0$ is the set of all horizontal subspaces~$H$ of
 $T_\lambda P^0$ such that their structure functions satisfy the chosen normalization condition $\mathcal{N}$.
 Obviously, the f\/ibers of $(P^0)^{(1)}$ are not empty, and if two horizontal subspaces $H$, $\tilde H$ belong to the f\/iber, then  $f_{H\widetilde H}\in \ker\, \partial$.
The subspace $\mg^1$ of ${\rm Hom}(\mg^{-1},\mg^{0})$ def\/ined by
 \begin{gather*}
 \mg^1\stackrel{\text{def}}{=} \ker \partial.
 \end{gather*}
 is called \emph{the first algebraic prolongation of $\mg^0\subset\mathfrak {gl}(\mg^{-1})$}. Note that it is absolutely not
 important that $\mg^0$ be a subalgebra of  $\mathfrak {gl}(\mg^{-1})$: the f\/irst algebraic prolongation can be def\/ined for a subspace of
 $\mathfrak {gl}(\mg^{-1})$ (see the further generalization below).

 If $\mg^1 =0$ then the choice of the ``normalization conditions'' $\mathcal N$ determines an Ehresmann
connection  on $P^0$ and $P^0$ is endowed with a canonical frame. As an example consider a~Riemannian structure. In this case $\mg^0=\mathfrak{so}(n)$, where
$n=\dim \mg^{-1}$, and
it is easy to show that $\mg^1=0$. Moreover, $\dim {\rm Hom}(\mg^{-1}\wedge\mg^{-1},\mg^{-1})=
\dim {\rm Hom}(\mg^{-1},\mg^{0})=\frac{n^2(n-1)}{2}$.
Hence, $\text{Im}\,\partial={\rm Hom}(\mg^{-1}\wedge\mg^{-1},\mg^{-1})$ and the complement subspace $\mathcal N$ must be equal to $0$.
So, in this case one gets the canonical Ehresmann connection with zero structure function (torsion), which is nothing but the
Levi-Civita connection.

If $\mg^1\neq 0$, we continue the prolongation procedure by induction.
 Given a linear space $W$ denote by ${\rm Id}_W$ the identity map on $W$. The bundle
 $(P^0)^{(1)}$ is a frame bundle with the abelian structure group $G^1$ of all maps $A\in {\rm GL}(\mg^{-1}\oplus\mg^{0})$
 such that
\begin{gather}
A|_{\mg^{-1}}={\rm Id}_{\mg^{-1}}+T, \nonumber\\
A|_{\mg^0}={\rm Id}_{\mg^0},\label{str}
\end{gather}
where $T\in \mg^1$. The right action $R_A$ of $A\in G_1$ on a f\/iber of $(P^0)^{(1)}$  is def\/ined by the following rule:
$R_A (\vf) =\vf\circ A$. Observe that $\mg^1$ is isomorphic to the Lie algebra of $G^1$.

Set $P^1=(P^0)^{(1)}$.
The second prolongation $P^2$ of $P^0$ is by def\/inition the f\/irst prolongation of the frame bundle $P^1$,
$P^2\stackrel{\text{def}}{=}(P^1)^{(1)}$ and so on by induction: the $i$-th prolongation $P^{i}$ is  the f\/irst prolongation
of the frame bundle $P^{i-1}$.

Let us describe the structure group $G^i$ of the frame bundle $P^i$ over $P^{i-1}$ in more detail.
For this one can def\/ine the
Spencer operator and the f\/irst algebraic prolongation also for a~subspace~$W$ of ${\rm Hom} (\mg^{-1},V)$, where $V$ is a linear space, which does not
necessary coincide with $\mg^{-1}$ as before. In this case the Spencer operator is the operator from ${\rm Hom} (\mg^{-1}, W)$ to
${\rm Hom}(\mg^{-1}\wedge \mg^{-1},
V)$, def\/ined by the same formulas, as in~\eqref{Spencer}. The f\/irst prolongation~$W^{(1)}$  of~$W$
is the kernel of the
Spencer operator. Note that by def\/inition $\mg^{1}=(\mg^0)^{(1)}$. Then the $i$-th prolongation~$\mg^i$ of~$\mg^0$ is
def\/ined by the following recursive formula: $\mg^i=(\mg^{i-1})^{(1)}$. Note that $\mg^i\subset {\rm Hom} (\mg^{-1},\mg^{i-1})$.
By \eqref{str} and the def\/inition of the Spencer operator the bundle
 $P^i$ is a~frame bundle with the abelian structure group $G^i$ of all maps
 $\displaystyle{A\in {\rm GL}\bigl(\bigoplus_{p=-1}^{i-1}\mg^{p}\bigr)}$ such that
\begin{gather*}
 A|_{\mg^-1}={\rm Id}_{\mg^{-1}}+T,\nonumber \\
 A|_{\bigoplus_{p=0}^{i-1}\mg^{p}}={\rm Id}_{\bigoplus_{p=0}^{i-1}\mg^{p}},
\end{gather*}
where $T\in \mg^i$. In particular, if $\mg^{l+1}=0$ for some $l\geq 0$, then the bundle $P^{l}$ is endowed with the canonical frame
and we are done.

\section{Tanaka's f\/irst prolongation}\label{section3}

Now consider the general case. As before $P^0$ is a structure of constant type $(\mathfrak m, \mg^0)$.
Let $\Pi_0:P^0\to M$ be the canonical projection. The f\/iltration $\{D^i\}_{i<0}$ of $TM$ induces a f\/iltration $\{D^i_0\}_{i\leq 0}$ of $T P^0$ as follows:
\begin{gather*}
 D^0_0=\ker (\Pi_0)_*,\nonumber\\
  D^i_0(\lambda)=\bigl\{v\in T_\lambda P^0: (\Pi_0)_*v\in D^i\bigl(\Pi_0(\lambda)\bigr)\bigr\}\qquad \forall\,  i<0.
\end{gather*}
We also set $D^i_0=0$ for all $i>0$.
Note that $D^0_0(\lambda)$ is the tangent space at $\lambda$ to the f\/iber of $P^0$ and therefore can be identif\/ied with $\mg^0$. Denote by $I_\lambda:\mg^0\to D^0_0(\lambda)$ the identifying isomorphism.

Fix a point $\lambda\in P^0$ and let  $\pi_0^i:D_0^i(\lambda)/ D_0^{i+2}(\lambda)\to D_0^i(\lambda)/ D_0^{i+1}(\lambda)$ be the canonical projection to the factor space.
Note that $\Pi_{0_*}$ induces an isomorphism between  the space $D_0^i(\lambda)/ D_0^{i+1}(\lambda)$ and the space $D^i(\Pi_0(\lambda))/ D^{i+1}(\Pi_0(\lambda))$ for any $i<0$. We denote this isomorphism by $\Pi_0^i$.
The f\/iber of the bundle $P^0$ over a point $x\in M$ is a subset of the set of all maps
\begin{gather*}\vf\in \bigoplus_{i<0} \text {Hom}\bigl(\mg^i,
D^i(x)/ D^{i+1}(x)\bigr),
\end{gather*}
 which are isomorphisms of the graded Lie algebras $\mathfrak{m}=\displaystyle{\bigoplus_{i<0} \mg^i}$ and
$\displaystyle{\bigoplus_{i<0} D^i(x)/D^{i+1}(x)}$.  We are going to construct a new bundle $P^1$
over the bundle $P^0$ such that the f\/iber of $P^1$ over a point $\lambda=(x,\vf)\in P^0$ will be a certain subset of
the set of all maps
\begin{gather*}\hat\vf\in \bigoplus_{i\leq 0} \text{Hom}\bigl(\mg^i,
D_0^i(\lambda)/ D_0^{i+2}(\lambda)\bigr)
\end{gather*}
 such that
\begin{gather}
 \vf|_{\mg^i}=\Pi_0^i\circ\pi_0^i\circ\hat\vf|_{\mg^i} \qquad \forall \, i<0,\nonumber\\
 \hat\vf|_{\mg^0}=I_\lambda.\label{hat}
\end{gather}
For this f\/ix again a point $\lambda=(x,\vf)\in P^0$. For any $i<0$ choose a subspace $H^i\subset D_0^i(\lambda)/ D_0^{i+2}(\lambda)$, which is a complement of
$D_0^{i+1}(\lambda)/D_0^{i+2}(\lambda)$ to $D_0^i(\lambda)/ D_0^{i+2}(\lambda)$:
\begin{gather}
\label{H1}
D_0^i(\lambda)/ D_0^{i+2}(\lambda)=D_0^{i+1}(\lambda)/ D_0^{i+2}(\lambda)\oplus H^i.
\end{gather}
Then the map $\Pi^i_0\circ\pi^i_0|_{H^i}$ def\/ines an isomorphism between $H^i$ and $D^i\bigl(\Pi_0(\lambda)\bigr)/ D^{i+1}\bigl(\Pi_0(\lambda)\bigr)$.  So, once a tuple of subspaces $\mathcal H=\{H^i\}_{i<0}$ is chosen,
one can def\/ine a map
\begin{gather*}
\vf^{\mathcal H}\in \displaystyle{\bigoplus_{i\leq 0} \text{Hom}\bigl(\mg^i,
D_0^i(\lambda)/ D_0^{i+2}(\lambda)\bigr)}
\end{gather*}
 as follows
\begin{gather*}
\vf^{\mathcal H}|_{\mg^i}=
\begin {cases}\big(\Pi^i_0\circ\pi_0^i|_{H^i}\big)^{-1}\circ\vf|_{\mg^i} & \text{ if } i<0,\\
I_\lambda & \text { if } i=0.
\end{cases}
\end{gather*}
Clearly $\hat\vf=\vf^{\mathcal H}$ satisf\/ies  \eqref{hat}. Tuples of subspaces $\mathcal H=\{H^i\}_{i<0}$ satisfying \eqref{H1}
play here the same role as horizontal subspaces in the prolongation of the usual $G$-structures. Can we choose a tuple
$\{H^i\}_{i<0}$ in a canonical way? For this, by analogy with the prolongation of $G$-structure, we introduce a ``partial soldering form'' of the bundle $P^0$ and the structure function of a~tuple~$\mathcal H$. The \emph{soldering form} of $P^0$ is a tuple
$\Omega_0=\{\omega_0^i\}_{i<0}$, where $\omega^i_0$ is a $\mg^i$-valued linear form on $D_0^i(\lambda)$
def\/ined by
\begin{gather*}
\omega_0^i(Y)=\vf^{-1}\bigl(\bigl((\Pi_0)_*(Y)\bigr)_i\bigr),
\end{gather*}
where $\bigl((\Pi_0)_*(Y)\bigr)_i$ is the equivalence class of $(\Pi_0)_*(Y)$ in $D^i(x)/D^{i+1}(x)$.
Observe that $D_0^{i+1}(\lambda)$ $=\ker \omega_0^i$. Thus the form $\omega_0^i$ induces the $\mg^i$-valued form $\bar \omega_0^i$ on $D_0^i(\lambda)
/D_0^{i+1}(\lambda)$.
The \emph {structure function $C_{\mathcal H}^0$ of the tuple $\mathcal H=\{H^i\}_{i<0}$} is the element of the space
\begin{gather}
\label{A0}
\mathcal A_0=
\left(\bigoplus_{i=-\mu}^{-2} {\rm Hom}(\mg^{-1}\otimes\mg^i,\mg^{i})\right)
\oplus {\rm Hom}\big(\mg^{-1}\wedge\mg^{-1},\mg^{-1}\big)
\end{gather}
 def\/ined as follows.
Let $\text{pr}_i^{\mathcal H}$ be the projection of $D_0^i(\lambda)/ D_0^{i+2}(\lambda)$ to $D_0^{i+1}(\lambda)/ D_0^{i+2}(\lambda)$
parallel to $H^i$ (or corresponding to the splitting \eqref{H1}). Given vectors $v_1\in \mg^{-1}$ and $v_2\in\mg^{i}$, take two vector f\/ields $Y_1$ and $Y_2$ in a neighborhood of $\lambda$ in $P^0$ such that $Y_1$ is a section of $D_0^{-1}$, $Y_2$ is a~section of $D_0^i$, and
\begin{gather}
 \omega_0^{-1}(Y_1)\equiv v_1,\qquad \omega_0^i(Y_2)\equiv v_2,\nonumber
\\
  Y_1(\lambda)=\vf^{\mathcal H}(v_1),\qquad Y_2(\lambda)\equiv \vf^{\mathcal H}(v_2)\,\,{\rm mod}\, D_0^{i+2}(\lambda).\label{Y1Y2}
\end{gather}
Then set
\begin{gather}
\label{structT1}
C_{\mathcal H}^0(v_1,v_2)\stackrel{\text{def}}{=}\bar\omega_0^i\bigl({\rm pr}_{i-1}^{\mathcal H}\bigl([Y_1,Y_2](\lambda)\bigr)\bigr).
\end{gather}
In the above formula we take the equivalence class of the vector $[Y_1, Y_2](\lambda)$ in $D_0^{i-1}(\lambda)/D_0^{i+1}(\lambda)$ and then apply ${\rm pr}_{i-1}^{\mathcal H}$.

One must show that $C_{\mathcal H}^0(v_1,v_2)$ does not depend on the choice of vector f\/ields~$Y_1$ and~$Y_2$, satisfying~\eqref{Y1Y2}.
Indeed, assume that $\widetilde Y_1$ and $\widetilde Y_2$ are another pair of vector f\/ields in a neighborhood of $\lambda$ in~$P^0$ such that $\widetilde Y_1$ is a section of~$D_0^{-1}$, $\widetilde Y_2$ is a section of $D_0^i$, and they satisfy~\eqref{Y1Y2}
with~$Y_1$,~$Y_2$ replaced by $\widetilde Y_1$, $\widetilde Y_1$. Then
\begin{gather}
\label{zz}
\widetilde Y_1=Y_1+Z_1, \qquad \widetilde Y_2=Y_2+Z_2,
\end{gather}
where $Z_1$ is a section of the distribution $D_0^0$ such that $Z_1(\lambda)=0$ and $Z_2$ is a section of the distribution $D_0^{i+1}$ such that $Z_2(\lambda)\in D_0^{i+2}(\lambda)$.
 It follows that $[Y_1, Z_2](\lambda) \in D_0^{i+1}(\lambda)$ and $[Y_2, Z_1](\lambda)\in D_0^{i+1}(\lambda)$. This together with the fact that $[Z_1,Z_2]$ is a section of $D_0^{i+1}$ imply that
\begin{gather*}
[\widetilde Y_1,\widetilde Y_2](\lambda)\equiv [Y_1,Y_2]\ \
\text{mod}\ D_0^{i+1}(\lambda).
 \end{gather*}
From
\eqref{structT1} we see that the structure function is independent
of the choice of vector f\/ields~$Y_1$ and~$Y_2$.


We now take another tuple $\widetilde {\mathcal H}=\{\widetilde H^i\}_{i<0}$ such that
\begin{gather}
\label{tildeH1}
D_0^i(\lambda)/ D_0^{i+2}(\lambda)=D_0^{i+1}(\lambda)/ D_0^{i+2}(\lambda)\oplus \widetilde H^i
\end{gather}
and consider how the structure functions $C_{\mathcal H}^1$ and $C_{\widetilde{\mathcal H}}^1$ are related.
By construction, for any vector $v\in\mg ^i$ the vector $\vf^{\widetilde{\mathcal H}}(v)-\vf^{\mathcal H}(v)$ belongs to
$D_0^{i+1}(\lambda)/ D_0^{i+2}(\lambda)$.
Let
\begin{gather*}
f_{\mathcal H\widetilde {\mathcal H}}(v)\stackrel{\text{def}}{=} \begin{cases}
\bar\omega_0^{i+1}\big(\vf^{\widetilde{\mathcal H}}(v)-\vf^{\mathcal H}(v)\big) & \text{ if } v\in \mg^i \text{ with } i<-1,\\
I_\lambda^{-1}\big(\vf^{\widetilde{\mathcal H}}(v)-\vf^{\mathcal H}(v)\big) & \text{ if } v\in \mg^{-1}.
\end{cases}
\end{gather*}
Then $f_{\mathcal H\widetilde {\mathcal H}}\in \displaystyle{\bigoplus_{i<0}{\rm Hom}(\mg^i,\mg^{i+1})}$.
Conversely, it is clear that for any $f\in\displaystyle{\bigoplus_{i<0}{\rm Hom}(\mg^i,\mg^{i+1})}$ there exists a tuple $\widetilde{\mathcal H}=\{\widetilde H^i\}_{i<0}$, satisfying \eqref{tildeH1}, such that  $f=f_{\mathcal H \widetilde {\mathcal H}}$.

Further, let $\mathcal A_0$ be as in \eqref{A0} and def\/ine
a map
\begin{gather*}
\partial_0:\displaystyle{\bigoplus_{i<0}{\rm Hom}(\mg^i,\mg^{i+1})}\rightarrow \mathcal A_0
\end{gather*}
  by
\begin{gather*}
\partial_0 f(v_1,v_2)=[f(v_1),v_2]+[v_1,f(v_2)]-f([v_1,v_2]),
\end{gather*}
where the brackets $[\, \,, \,]$ are as in the Lie algebra $\mathfrak m\oplus\mg^0$.
The map $\partial_0$ coincides with the Spencer operator \eqref{Spencer} in the case of $G$-structures. Therefore it is called the \emph{generalized Spencer operator for the first prolongation}.
\begin{proposition}
\label{prop1}
The following identity holds
\begin{gather}
\label{structrans}
C_{\widetilde {\mathcal H}}^0=C_{\mathcal H}^0+\partial_0f_{\mathcal H \widetilde{\mathcal H}}.
\end{gather}
\end{proposition}

\begin{proof}
Fix vectors $v_1\in \mg^{-1}$ and $v_2\in\mg^{i}$ and let $Y_1$ and $Y_2$ be two vector f\/ields in a neighborhood of $\lambda$ satisfying~\eqref{Y1Y2}. Take two vector f\/ields $\widetilde Y_1$ and $\widetilde Y_2$ in a neighborhood of $\lambda$ in~$P^0$ such that~$\widetilde Y_1$ is a section of $D_0^{-1}$, $\widetilde Y_2$ a section of $D_0^i$, and
\begin{gather*}
 \omega_0^{-1}(\widetilde Y_1)\equiv v_1,\qquad  \omega_0^i(\widetilde Y_2)\equiv v_2,
\nonumber\\
  \widetilde Y_1(\lambda)=\vf^{\widetilde{\mathcal H}}(v_1),\qquad \widetilde Y_2(\lambda)\equiv \vf^{\widetilde{\mathcal H}}(v_2)\,\,{\rm mod}\,\, D_0^{i+2}(\lambda).
\end{gather*}
Further, assume that vector f\/ields $Z_1$ and $Z_2$ are def\/ined as in \eqref{zz}. Then $Z_1$ is a section of $D_0^0$ and $Z_2$ is a section of $D_0^{i+1}$ such that
\begin{gather}
 Z_1(\lambda)=I_\lambda\bigl(f_{\mathcal H\widetilde {\mathcal H}}(v_1)\bigr), \label{zz11}\\
 f_{\mathcal H\widetilde {\mathcal H}}(v_2)=\begin{cases}
\bar\omega_0^{i+1}\left(Z_2(\lambda)\right) & \text{ if } v\in \mg^i,
\  i<-1,\\
I_\lambda^{-1}\left(Z_2(\lambda)\right) & \text{ if } v\in \mg^{-1}.
\end{cases}
\label{zz12}
\end{gather}
Hence $[Z_1, Y_2]$ and $[Y_1, Z_2]$ are sections of $D_0^i$, while $[Z_1, Z_2]$ is a section of $D_0^{i+1}$. This implies that
\begin{gather}
\label{inter1}
\bar\omega_0^i\Bigl({\rm pr}_{i-1}^{\widetilde {\mathcal H}}\bigl([\widetilde Y_1,\widetilde Y_2](\lambda)\bigr)\Bigr)=
\bar\omega_0^i\Bigl({\rm pr}_{i-1}^{\widetilde{\mathcal H}}\bigl([Y_1, Y_2](\lambda)\bigr)\Bigr)+\bar\omega_0^i\bigl([Z_1, Y_2]\bigr) + \bar\omega_0^i\bigl([Y_1, Z_2]\bigr).
\end{gather}
Further, directly from the def\/initions of $f_{\mathcal H\widetilde {\mathcal H}}$, ${\rm pr}_{i-1}^{\mathcal H}$, and ${\rm pr}_{i-1}^{\widetilde {\mathcal H}}$ it follows that
\begin{gather}
\label{inter2}
\bar\omega_0^i\bigl({\rm pr}_{i-1}^{\widetilde {\mathcal H}}(w)\bigr)=\bar\omega_0^i\bigl({\rm pr}_{i-1}^{\mathcal H}(w)\bigr)-f_{\mathcal H\widetilde {\mathcal H}}\bigl(\bar\omega_0^{i-1}(w)\bigr) \qquad \forall \, w\in D_0^{i-1}(\lambda)/D_0^{i+1}(\lambda).
\end{gather}
Besides, from the def\/inition of the soldering form, the fact that $\vf$ is an isomorphism of the Lie algebras $\mathfrak m$ and $\mathfrak m(x)=
\displaystyle{\bigoplus_{i<0}D^i(x)/D^{i+1}(x)}$ , and relations \eqref{zz12} for $i<-1$ it follows that
\begin{gather}
\bar\omega_0^i\bigl([Y_1, Z_2]\bigr)=[v_1, f_{\mathcal H\widetilde {\mathcal H}}(v_2)]\qquad \forall\, i<-1.
\label{inter3}
\end{gather}
Taking into account  \eqref{Y1Y2} we get
\begin{gather}
\bar\omega_0^{i-1}([Y_1, Y_2])=[v_1,v_2].
\label{inter4}
\end{gather}
Finally, from \eqref{zz11} and \eqref{zz12} for $i=-1$, and the def\/inition of the action of $G^0$ on $P^0$ it follows that
identity \eqref{inter3} holds also for $i=-1$, and that
\begin{gather}
\bar\omega_0^i\bigl([Z_1, Y_2]\bigr)=[f_{\mathcal H\widetilde {\mathcal H}}(v_1), v_2].
\label{inter5}
\end{gather}
Substituting \eqref{inter2}--\eqref{inter5} into \eqref{inter1} we get \eqref{structrans}.
\end{proof}

Now we proceed as in the case of $G$-structures.  Fix a subspace
\begin{gather*}
\mathcal N_0\subset \mathcal A_0
 \end{gather*}
  which is complementary to  $\text{Im}\, \partial_0$,
\begin{gather}
\label{normsplit}
\mathcal A_0= \text{Im} \,\partial_0\oplus\mathcal N_0.
 \end{gather}
 As for $G$-structures, the subspace $\mathcal N_0$ def\/ines the normalization conditions for the f\/irst prolongation. Then from the splitting \eqref{normsplit} it follows trivially  that there exists a tuple $\mathcal H=\{H^i\}_{i<0}$ such that
 \begin{gather}
 \label{norm0}
 C_{\mathcal H}^0\in \mathcal N_0.
 \end{gather}
  A tuple $\widetilde {\mathcal H}=\{\widetilde H^i\}_{i<0}$ satisf\/ies $C_{\widetilde{\mathcal H}}^0\in \mathcal N_0$ if and only if $f_{\mathcal H\widetilde{\mathcal H}}\in \ker\, \partial_0$. In particular if $\ker \,\partial_0 =0$ then the tuple $\mathcal H$ is f\/ixed uniquely by condition \eqref{norm0}.
 Let
 \begin{gather*}
 \mg^1\stackrel{\text{def}}{=} \ker \partial_0.
 \end{gather*}
 The space $\mg^1$ is called the \emph{first algebraic prolongation} of the algebra $\mathfrak m\oplus \mg^0$.
 Here we consider~$\mg^1$ as an abelian Lie algebra. Note that the fact that the symbol $\mathfrak m$ is fundamental (that is, $\mg^{-1}$~generates the whole $\mathfrak m$) implies that
 \begin{gather*}
\mg^1=\left\{f\in \bigoplus_{i<0}{\rm Hom}(\mg^i,\mg^{i+1}): f ([v_1,v_2])=[f (v_1),v_2]+[v_1, f( v_2)]\,\,\forall\, v_1, v_2 \in \mathfrak m\right\}.
\end{gather*}
 The \emph{first $($geometric$)$  prolongation} of the bundle $P^0$ is the bundle $P^1$ over $P^0$ def\/ined by
 \begin{gather*}
 P^1=\big\{(\lambda, \mathcal H):\lambda\in P^0, C_{\mathcal H}^0\in \mathcal N_0\big\}.
 \end{gather*}
Equivalently,
\begin{gather*}
 P^1=\big\{(\lambda,\vf^{\mathcal H}):\lambda\in P^0, C_{\mathcal H}^0\in \mathcal N_0\big\}.
 \end{gather*}
It is a principal bundle with the abelian structure group $G^1$ of all maps $A\in\! \displaystyle{\bigoplus_{i<1}{\rm Hom}(\mg^i,\mg^{i}\oplus\mg^{i+1})}$ such that
\begin{gather*}
 A|_{\mg^i}={\rm Id}_{\mg^i}+T_i, \qquad i<0,\nonumber\\
 A|_{\mg^0}={\rm Id}_{\mg^0},
\end{gather*}
where $T_i\in {\rm Hom}(\mg^i,\mg^{i+1})$ and $(T_{-\mu},\ldots,T_{-1})\in \mg^1$. The right action $R_A^1$ of $A\in G^1$ on a f\/iber of $P^1$  is def\/ined by $R_A^1 (\vf^{\mathcal H}) =\vf^{\mathcal H}\circ A$. Note that $G^1$ is an abelian group of dimension equal to $\dim \mg^1$.

\section{Higher order Tanaka's prolongations}\label{section4}

More generally,  def\/ine the $k$-th algebraic prolongation $\mg^k$ of the algebra $\mathfrak m\oplus \mg^0$ by induction for any $k\in\mathbb N$. Assume that spaces $\mg^l\subset \displaystyle{\bigoplus_{i<0}{\rm Hom}(\mg^i,\mg^{i+l})}$ are def\/ined  for
all $0<l<k$. Set
\begin{gather}
\label{br2}
[f,v]=-[v, f]=f(v) \qquad \forall \, f\in \mg^l, \  0\leq l<k,  \ \text{and} \ v\in\mathfrak m.
\end{gather}
 Then let
\begin{gather}
\label{mgk}
\mg^k\stackrel{\text{def}}{=}\left\{f\in \bigoplus_{i<0}{\rm Hom}(\mg^i,\mg^{i+k}): f ([v_1,v_2])=[f (v_1),v_2]+[v_1, f( v_2)]\,\,\forall\, v_1, v_2 \in \mathfrak m\right\}.
\end{gather}
Directly from this def\/inition and the fact that $\mathfrak m$ is fundamental (that is, it is generated by $\mg^{-1}$) it follows that if
$f\in \mg^k$ satisf\/ies $f|_{\mg^{-1}}=0$, then $f=0$.
The space $\bigoplus_{i\in Z} \mg^i$ can be naturally endowed with the structure of a graded Lie algebra.
The brackets of two elements from $\mathfrak m$ are as in $\mathfrak m$. The brackets of an element with non-negative weight and an element from $\mathfrak m$ are already def\/ined by  \eqref{br2}.
It only remains to def\/ine  the brackets $[f_1,f_2]$ for $f_1\in\mg^k$, $f_2\in \mg^l$ with $k,l\geq 0$.
The def\/inition is inductive with respect to $k$ and $l$: if $k=l=0$ then the bracket $[f_1,f_2]$ is as in~$\mg^0$. Assume that $[f_1,f_2]$
is def\/ined for all $f_1\in\mg^k$, $f_2\in \mg^l$ such that a pair  $(k,l)$ belongs to  the set
\begin{gather*}
\{(k,l):0 \leq k\leq \bar k, 0\leq l\leq \bar l\}\backslash \{(\bar k,\bar l)\}.
\end{gather*}
Then def\/ine $[f_1, f_2]$  for $f_1\in\mg^{\bar k}$, $f_2\in \mg^{\bar l}$ to be the element of $\displaystyle{\bigoplus_{i<0}{\rm Hom}(\mg^i,\mg^{i+\bar k+\bar l})}$ given by
\begin{gather}
\label{posbr}
[f_1,f_2]v\stackrel{\text{def}}{=} [f_1(v), f_2]+[f_1,f_2(v)]\qquad \forall \, v\in\mathfrak m.
\end{gather}
It is easy to see that $[f_1,f_2]\in \mg^{k+l}$ and that $\bigoplus_{i\in Z} \mg^i$ with bracket product def\/ined as above is a~graded Lie algebra. As a matter of fact \cite[\S~5]{tan} this graded Lie algebra satisf\/ies Properties~1--3 from Subsection~\ref{section1.4}. That is it is a realization of the algebraic universal prolongation $\mg(\mathfrak m, \mg^0)$ of the algebra $\mathfrak m\oplus\mg^0$.

Now we are ready to construct the higher order geometric prolongations of the bundle $P^0$ by induction.  Assume that all  $l$-th order prolongations $P^l$ are constructed for $0\leq l\leq k$. We also set $P^{-1}=M$. We will not specify what the bundles $P^l$ are exactly. As in the case of the f\/irst prolongation $P^1$,  their construction depends on the choice of normalization conditions on each step. But we will point out those properties of these bundles that we need in order to construct the $(k+1)$-st order prolongation $P^{k+1}$. Here are these properties:
\begin{enumerate}\itemsep=0pt
\item $P^l$ is a principal bundle over $P^{l-1}$ with an abelian structure group $G^l$ of dimension  equal to $\dim \mg^l$ and with  the canonical projection $\Pi_l$.
\item
The tangent bundle $T P^l$ is endowed with the f\/iltration $\{D^i_l\}$ as follows: For $l=-1$ it coincides with the initial f\/iltration  $\{D^{i}\}_{i<0}$ and  for $l\geq 0$
we get by induction
\begin{gather*}
 D^l_l=\ker (\Pi_l)_*,\nonumber\\
  D^i_l(\lambda_l)=\bigl\{v\in T_\lambda P^l: (\Pi_l)_*v\in D_{l-1}^i\bigl(\Pi_l(\lambda_l)\bigr)\bigr\}  \qquad \forall\, i<l.
\end{gather*}
The subspaces $D^l_l(\lambda_l)$, as the tangent spaces to the f\/ibers of $P^l$ , are canonically identif\/ied with $\mg^l$. Denote by $I_{\lambda_l}:\mg^l\rightarrow D^l_l(\lambda_l)$ the identifying isomorphism.
\item The f\/iber of $P^l$, $0\leq l\leq k$, over a point $\lambda_{l-1}\in P^{l-1}$ will be a certain subset of
the set of all maps from
\begin{gather*}
\bigoplus_{i< l} \text{Hom}\bigl(\mg^i,
D_{l-1}^i(\lambda_{l-1})/ D_{l-1}^{i+l+1}(\lambda_{l-1})\bigr).
\end{gather*}
  If $l>0$ and $\lambda_l=(\lambda_{l-1}, \vf_l) \in P^l$, then $\vf_l|_{\mg^{l-1}}$ coincides with the identif\/ication of $\mg^{l-1}$ with $D^{l-1}_{l-1}(\lambda_{l-1})$ and the restrictions $\vf_l|_{\mg^i}$ with $i\geq 0$ are the same  for all $\lambda_l$ from the same f\/iber.

\item Assume that $0<l\leq k$, $\lambda_{l-1}=(\lambda_{l-2}, \vf_{l-1})\in P^{l-1}$ and $\lambda_l=(\lambda_{l-1},\vf_l)\in P^l$.
    The maps $\vf_{l-1}$ and $\vf_l$ are related as follows: if
    \begin{gather}
    \label{pli}
    \pi_l^i: \ D_l^i(\lambda_{l})/ D_l^{i+l+2}(\lambda_{l})\to D_l^i(\lambda_{l})/ D_l^{i+l+1}(\lambda_{l})
     \end{gather}
     are the canonical projections to a factor space and
     \begin{gather}
     \label{Pli}
     \Pi_l^i: \ D^i_l(\lambda_l)/D^{i+l+1}_l(\lambda_l)\rightarrow D^i_{l-1}\bigl(\Pi_l(\lambda_l)\bigr)/D^{i+l+1}_{l-1}(\Pi_l(\lambda_l))
     \end{gather}
are the canonical maps induced by $(\Pi_l)_*$, then
\begin{gather*}
\forall \, i<l \qquad \vf_{l-1}|_{\mg^i}=\Pi_{l-1}^i\circ \pi_{l-1}^i\circ \vf_{l}|_{\mg^i}.
\end{gather*}
Note that the maps $\Pi_l^i$ are isomorphisms for $i<0$ and the maps $\pi_l^i$ are identities for $i\geq 0$ (we set $D_l^i=0$ for $i>l$).
\end{enumerate}

Now we are ready to construct the $(k+1)$-st order Tanaka geometric prolongation. Fix a point $\lambda_k \in P^k$ and assume that $\lambda_k=(\lambda_{k-1},\vf_k)$, where
\begin{gather*}
\vf_k\in \displaystyle{\bigoplus_{i< k} \text{Hom}\bigl(\mg^i,
D_{k-1}^i(\lambda_{k-1})/ D_{k-1}^{i+k+1}(\lambda_{k-1})\bigr)}.
\end{gather*}
 Let $\mathcal H_k=\{H_k^i\}_{i<k}$ be the tuple of spaces such that $H_k^i=\vf_k(\mg^i)$. Take a tuple  $\mathcal H_{k+1}=\{H_{k+1}^i\}_{i<k}$ of linear spaces such that
\begin{enumerate}\itemsep=0pt
\item for $i<0$ the space
$H_{k+1}^i$ is a complement
of $D_k^{i+k+1}(\lambda_k)/D_k^{i+k+2}(\lambda_k)$ in $(\Pi_k^i\circ\pi_k^i)^{-1} (H_k^i)\subset D_{k}^i(\lambda_{k})/ D_{k}^{i+k+2}(\lambda_{k})$,
\begin{gather}
\label{Hk-}
(\Pi_k^i\circ\pi_k^i)^{-1} (H_k^i)=D_k^{i+k+1}(\lambda_k)/D_k^{i+k+2}(\lambda_k)\oplus H_{k+1}^i;
\end{gather}

\item for $0\leq i<k$ the space $H_{k+1}^i$ is a complement of $D_k^k(\lambda_k)$ in  $(\Pi_k^i)^{-1} (H_k^i)$,
\begin{gather}
\label{Hk+}
(\Pi_k^i)^{-1} (H_k^i)=D_k^{k}(\lambda_k)\oplus H_{k+1}^i.
\end{gather}
\end{enumerate}
Here the maps $\pi_k^i$ and $\Pi_k^i$ are def\/ined as in \eqref{pli} and \eqref{Pli} with $l=k$.

Since $D_k^{i+k+1}(\lambda_k)/D_k^{i+k+2}(\lambda_k)=\ker \pi_k^i$ and $\Pi_k^i$ is an isomorphism for $i<0$,  the map $\Pi_k^i\circ\pi^i_k|_{H_{k+1}^i}$ def\/ines an isomorphism between $H_{k+1}^i$ and $H_k^i$ for $i<0$.
Additionally,  by \eqref{Hk+} the map $(\Pi_l)_*|_{H_{k+1}^i}$ def\/ines an isomorphism between $H_{k+1}^i$ and $H_k^i$ for $0\leq i<k$.
 So, once a tuple of subspaces $\mathcal H_{k+1}=\{H_{k+1}^i\}_{i<k}$, satisfying \eqref{Hk-} and \eqref{Hk+}, is chosen,
one can def\/ine a map
\begin{gather*}
\vf^{\mathcal H_{k+1}}\in \bigoplus_{i\leq k} \text{Hom}\bigl(\mg^i,
D_k^i(\lambda_k)/ D_k^{i+k+2}(\lambda_k)\bigr)
\end{gather*}
 as follows
\begin{gather*}
\vf^{\mathcal H_{k+1}}|_{\mg^i}=
\begin{cases}
\big(\Pi_k^i\circ\pi^i_k|_{H_{k+1}^i}\big)^{-1}\circ\vf_k|_{\mg^i} &\text { if } i<0,\\
\bigl((\Pi_l)_*|_{H_{k+1}^i}\bigr)^{-1}\circ\vf_k|_{\mg^i} &\text{ if }0\leq i<k,\\
I_{\lambda_k}& \text{ if } i=k.
\end{cases}
\end{gather*}
Can we choose a tuple or a subset of tuples
 $\mathcal H_k$ in a canonical way? To answer this question, by analogy with Sections~\ref{section2} and~\ref{section3}, we introduce a ``partial soldering form'' of the bundle $P^k$ and the structure function of a tuple $\mathcal H_{k+1}$.
The \emph{soldering form} of $P^k$ is a tuple
$\Omega_k=\{\omega_k^i\}_{i<k}$, where~$\omega_k^i$ is a $\mg^i$-valued linear form on $D_k^i(\lambda_k)$
def\/ined by
\begin{gather*}
\omega_k^i(Y)=\vf_k^{-1}\bigl(\bigl((\Pi_k)_*(Y)\bigr)_i\bigr).
\end{gather*}
Here $\bigl((\Pi_k)_*(Y)\bigr)_i$ is the equivalence class of $(\Pi_k)_*(Y)$ in $D_{k-1}^i(\lambda_{k-1})/D_{k-1}^{i+k+1}(\lambda_{k-1})$.
By construction it follows immediately that  $D_k^{i+1}(\lambda_k)=\ker \omega_k^i$. So, the form $\omega_k^i$ induces the $\mg^i$-valued form $\bar \omega_k^i$ on $D_{k}^i(\lambda_k)
/D_k^{i+1}(\lambda_k)$.

The \emph{ structure function $C_{\mathcal H_{k+1}}^k$ of a tuple $\mathcal H_{k+1}$} is the element of the space
\begin{gather}
{\mathcal A}_k=\left(\bigoplus_{i=-\mu}^{-2} {\rm Hom}\big(\mg^{-1}\otimes\mg^i,\mg^{i+k}\big)\right)
\oplus {\rm Hom}\big(\mg^{-1}\wedge\,\mg^{-1},\mg^{k-1}\big)\nonumber\\
\phantom{{\mathcal A}_k=}{} \oplus\left( \bigoplus_{i=0}^{k-1} {\rm Hom}\big(\mg^{-1}\otimes\mg^i,\mg^{k-1}\big)\right)\label{Ak}
\end{gather}
 def\/ined as follows:
Let  $\pi_l^{i,s}: D_l^i(\lambda_l)/D_l^{i+l+2}(\lambda_l)\rightarrow D_l^i(\lambda_l)/ D_l^{i+l+2-s}(\lambda_l)$ be the canonical projection to a factor space, where $-1\leq l\leq k$, $i\leq l$.
Here, as before, we assume that $D_l^i=0$ for $i>l$.
Note that the previously def\/ined $\pi_l^i$ coincides with $\pi_l^{i,1}$.
By construction, one has the following two relations
\begin{gather}
D_{k}^i(\lambda_{k})/D_{k}^{i+k+2}(\lambda_{k})=\left(\bigoplus_{s=0}^k\pi_k^{i+s,s} (H_{k+1}^{i+s})\right)\oplus
D_{k}^{i+k+1}(\lambda_{k})/D_{k}^{i+k+2}(\lambda_{k})
\quad \text{ if} \ \ i<0,
\label{longsplit1}
\\
 D_{k}^i(\lambda_{k})=\left(\bigoplus_{s=i}^{k-1}
H_{k+1}^{i}\right)\oplus D_{k}^{k}(\lambda_{k})\quad \text{if} \ \ 0\leq i<k.
\label{longsplit2}
\end{gather}
Let $\text{pr}_i^{\mathcal H_{k+1}}$ be the projection of $D_k^i(\lambda_k)/ D_k^{i+k+2}(\lambda_k)$ to $D_k^{i+k+1}(\lambda_k)/ D_k^{i+k+2}(\lambda_k)$
correspon\-ding to the splitting \eqref{longsplit1}
if $i<0$ or the projection of $D_k^i(\lambda_k)$ to
$H_{k+1}^{k-1}$
corresponding to the splitting \eqref{longsplit2} if $0\leq i<k$.
Given vectors $v_1\in \mg^{-1}$ and $v_2\in\mg^{i}$ take two vector f\/ields $Y_1$ and $Y_2$ in a neighborhood $U_k$ of $\lambda_k$ in $P^k$ such that for any $\tilde \lambda_k=(\tilde \lambda_{k-1},\tilde\vf_k)\in U_k$, where $\tilde \vf_k\in \displaystyle{\bigoplus_{i< k} \text{Hom}\bigl(\mg^i,
D_{k-1}^i(\lambda_{k-1})/ D_{k-1}^{i+k+1}(\lambda_{k-1})\bigr)}$, one has
\begin{gather}
\Pi_{k_{*}}Y_1(\tilde \lambda_k)= \tilde \vf_{k}(v_1), \qquad \Pi_{k_{*}}Y_2(\tilde\lambda_k)\equiv \tilde\vf_{k}(v_2) \quad {\rm mod}\,\, D_{k-1}^{i+k+1}(\lambda_{k-1}),\nonumber
\\
Y_1(\lambda)=\vf^{\mathcal H_{k+1}}(v_1),\qquad Y_2(\lambda)\equiv\vf^{\mathcal H_{k+1}}(v_2)\quad {\rm mod}\,\, D_k^{i+k+2}(\lambda).\label{Y1Y2k}
\end{gather}
Then set
\begin{gather}
\label{structTk}
C_{\mathcal H^{k+1}}^k(v_1,v_2)\stackrel{\text{def}}{=}
\begin{cases}
\bar\omega_k^{i+k}\bigl({\rm pr}_{i-1}^{\mathcal H_{k+1}}\bigl([Y_1,Y_2]\bigr)\bigr)
& \text{ if } i<0,\\
\omega_k^{k-1}\bigl({\rm pr}_{i-1}^{\mathcal H_{k+1}}\bigl([Y_1,Y_2]\bigr)\bigr)& \text{ if } 0\leq i<k.
\end{cases}
\end{gather}

As in the  case of the f\/irst prolongation, $C_{\mathcal H}^k(v_1,v_2)$ does not depend on the choice of vector f\/ields $Y_1$ and $Y_2$, satisfying \eqref{Y1Y2k}.
Indeed, assume that $\widetilde Y_1$ and $\widetilde Y_2$ is another pair of vector f\/ields in a neighborhood of $\lambda_k$ in $P^k$ such that $\widetilde Y_1$ is a section of $D_k^{-1}$, $\widetilde Y_2$ is a section of $D_k^i$, and they satisfy \eqref{Y1Y2k}
with  $Y_1$, $Y_2$ replaced by $\widetilde Y_1$, $\widetilde Y_1$.
Then
\begin{gather*}
\widetilde Y_1=Y_1+Z_1, \qquad \widetilde Y_2=Y_2+Z_2,
\end{gather*}
where $Z_1$ is a section of the distribution $D_k^k$ such that $Z_1(\lambda_k)=0 $ and $Z_2$ is a section of the distribution $D_k^
{\min\{i+k+1,k\}}$ such that $Z_2(\lambda_k)\in D_k^
{\min\{i+k+1,k\}+1}
(\lambda_k)$.
 Then $[Y_1, Z_2](\lambda_k) \in D_k^
 {\min\{i+k+1,k\}}
 (\lambda_k)$ and $[Y_2, Z_1](\lambda)\in D_k^
 {\min\{i+k+1,k\}}
 (\lambda_k)$. This together with the fact that $[Z_1,Z_2]$ is a section of $D_k^
 {\min\{i+k+1,k\}+1}
 $ implies that
\begin{gather*}
[\widetilde Y_1,\widetilde Y_2](\lambda)\equiv [Y_1,Y_2]\quad
\text{mod}\,\,D_k^
{\min\{i+k+1,k\}}
(\lambda).
 \end{gather*}
From \eqref{structTk} it follows that the structure function is independent
of the choice of vector f\/ields~$Y_1$ and~$Y_2$.

Now take another tuple $\widetilde {\mathcal H}_{k+1}=\{\widetilde H_{k+1}^i\}_{i<k}$ such that
\begin{enumerate}
\item for $i<0$ the space
$\widetilde H_{k+1}^i$ is a complement
of $D_k^{i+k+1}(\lambda_k)/D_k^{i+k+2}(\lambda_k)$ in $(\Pi_k^i\circ\pi_k^i)^{-1} (H_k^i)\subset D_{k}^i(\lambda_{k})/ D_{k}^{i+k+2}(\lambda_{k})$,
\begin{gather}
\label{Hk-t}
\big(\Pi_k^i\circ\pi_k^i\big)^{-1} (H_k^i)=D_k^{i+k+1}(\lambda_k)/D_k^{i+k+2}(\lambda_k)\oplus \widetilde H_{k+1}^i;
\end{gather}

\item  for $0\leq i<k$ the space $\widetilde H_{k+1}^i$ is a complement of $D_k^k(\lambda_k)$ in  $(\Pi_k^i)^{-1} (H_k^i)$,
\begin{gather}
\label{Hk+t}
\big(\Pi_k^i\big)^{-1} (H_k^i)=D_k^{k}(\lambda_k)\oplus \widetilde H_{k+1}^i.
\end{gather}
\end{enumerate}

How are the structure functions $C_{{\mathcal H}_{k+1}}^k$ and $C_{\widetilde{\mathcal H}_{k+1}}^k$ related?
By construction, for any vector $v\in\mg ^i$ the vector $\vf^{\widetilde{\mathcal H}_{k+1}}(v)-\vf^{\mathcal H_{k+1}}(v)$ belongs to
$D_k^{i+k+1}(\lambda_k)/ D_k^{i+k+2}(\lambda_k)$, for $i<0$, and to~$D_k^{k}(\lambda_k)$,  for $0\leq i<k$.
Let
\begin{gather*}
f_{\mathcal H_{k+1}\widetilde {\mathcal H}_{k+1}}(v)\stackrel{\text{def}}{=} \begin{cases}
\bar\omega_k^{i+k+1}\big(\vf^{\widetilde{\mathcal H}_{k+1}}(v)-\vf^{\mathcal H_{k+1}}(v)\big) & \text{ if } v\in \mg^i \text{ with } i<-1,\\
I_\lambda^{-1}\big(\vf^{\widetilde{\mathcal H}_{k+1}}(v)-\vf^{\mathcal H_{k+1}}(v)\big) & \text{ if } v\in \mg^{i} \text{ with } -1\leq i<k.
\end{cases}
\end{gather*}

Then
\begin{gather*}
f_{\mathcal H_{k+1}\widetilde {\mathcal H}_{k+1}}\in \displaystyle{\bigoplus_{i<0}{\rm Hom}(\mg^i,\mg^{i+k+1})}\oplus
\displaystyle{\bigoplus_{i=0}^{k-1}{\rm Hom}(\mg^i,\mg^k)}.
\end{gather*}
 In the opposite direction, it is clear that for any $f\in\displaystyle
 {\bigoplus_{i<0}{\rm Hom}(\mg^i,\mg^{i+k+1})}\oplus
\displaystyle{\bigoplus_{i=0}^{k-1}{\rm Hom}(\mg^i,\mg^k)}$, there exists a tuple $\widetilde{\mathcal H}_{k+1}=\{\widetilde H_{k+1}^i\}_{i<k}$ satisfying \eqref{Hk-t} and \eqref{Hk+t} and such that  $f=f_{\mathcal H_{k+1} \widetilde {\mathcal H}_{k+1}}$. Further, let $\mathcal A_k$ be as in \eqref{Ak} and def\/ine
a map
\begin{gather*}
\partial_k:\displaystyle{\bigoplus_{i<0}{\rm Hom}(\mg^i,\mg^{i+k+1})}\oplus
\displaystyle{\bigoplus_{i=0}^{k-1}{\rm Hom}(\mg^i,\mg^k)}\rightarrow \mathcal A_k
\end{gather*}
  by
\begin{gather}
\label{Spk}
\partial_k f(v_1,v_2)\\
\qquad{} =\begin{cases}[f(v_1),v_2]+[v_1,f(v_2)]-f([v_1,v_2])&\text{ if } v_1\in \mg^
{-1}, \ v_2\in \mg^i, \ i<0, \\
[v_1,f(v_2)] &\text{ if } v_1\in \mg^
{-1}, \ v_2\in \mg^i, \ 0\leq i
<k-1,
\end{cases}\nonumber
\end{gather}
where the brackets $[\, \,, \,]$ are as
  in the algebraic universal  prolongation $\mg(\mathfrak m,\mg^0)$.
For $k=0$ this def\/inition coincides with the def\/inition of the generalized Spencer operator for the f\/irst prolongation given in the previous section.

The reason for introducing the operator $\partial_k$ is that the following generalization of identity~\eqref{structrans} holds:
\begin{gather*}
C_{\widetilde {\mathcal H}_{k+1}}^k=C_{\mathcal H_{k+1}}^k+\partial_kf_{\mathcal H_{k+1} \widetilde{\mathcal H}_{k+1}}.
\end{gather*}

A verif\/ication of this identity for pairs $(v_1,v_2)$, where $v_1\in \mg^{-1}$ and
$v_2\in \mg^i$ with $i<0$, is completely analogous to the proof of Proposition~\ref{prop1}. For $i\geq 0$ one has to use the inductive assumption that the restrictions $\vf_l|_{\mg^i}$
are the same  for all $\lambda_l$ from the same f\/iber
(see item~3 from the list of properties satisf\/ied by $P^l$ in the beginning of this section) and the splitting \eqref{longsplit2}.

Now we proceed as in Sections~\ref{section2} and~\ref{section3}.  Fix a subspace
\begin{gather*}
\mathcal N_k\subset\mathcal A_k
 \end{gather*}
  which is complementary to  $\text{Im}\, \partial_k$,
\begin{gather}
\label{normsplitk}
\mathcal A_k= \text{Im} \,\partial_k\oplus\mathcal N_k.
 \end{gather}
 As above, the subspace $\mathcal N_k$ def\/ines the normalization conditions for the f\/irst prolongation. Then from the splitting~\eqref{normsplitk} it follows trivially  that there exists a tuple $\mathcal H_{k+1}=\{H_{k+1}^i\}_{i<k}$, satis\-fying~\eqref{Hk-} and~\eqref{Hk+}, such that
 \begin{gather*}
 C_{\mathcal H_{k+1}}^k\in \mathcal N_k
 \end{gather*}
  and
 $C_{\widetilde{\mathcal H}_{k+1}}^k\in \mathcal N_k$ for a tuple $\widetilde H_{k+1}=\{H_{k+1}^i\}_{i<k}$, satisfying \eqref{Hk-t} and \eqref{Hk+t}, if and only if $f_{\mathcal H_{k+1}\widetilde {\mathcal H}_{k+1}}\in \ker\, \partial_k$.
 Note also that
 \begin{gather}
 \label{kerk}
 f\in \ker\, \partial_k \,\Rightarrow\, f|_{\mg^i}=0  \qquad \forall \, 0\leq i\leq k-1.
 \end{gather}
 In other words,
 \begin{gather}
 \label{kerk1}
 \ker\, \partial_k\subset \displaystyle{\bigoplus_{i<0}{\rm Hom}\big(\mg^i,\mg^{i+k+1}\big)}.
 \end{gather}
 Indeed, if $f\in \ker\, \partial_k$, then by \eqref{Spk} for any $v_1\in\mg^{-1}$ and $v_2\in \mg^{i}$ with $0\leq i\leq k-1$ one has
 \begin{gather*}
 [v_1,f(v_2)]=-f(v_2)v_1=0.
 \end{gather*}
In other words, $f(v_2)|_{\mg^{-1}}=0$ (recall that $f(v_2)\in \mg^k\subset  \displaystyle{\bigoplus_{i<0}{\rm Hom}(\mg^i,\mg^{i+k})}$). Since $\mg^{-1}$ generates the whole symbol $\mathfrak m$ we see that $f(v_2)=0$ holds for any $v_2\in \mg^{i}$ with $0\leq i\leq k-1$. This proves that \eqref{kerk}.

 Further, comparing  \eqref{Spk} and \eqref{kerk1} with \eqref{mgk} and using again the fact that $\mg^{-1}$ generates the whole symbol $\mathfrak m$ we obtain
 \begin{gather*}
 \ker \partial_k=\mg^{k+1}.
 \end{gather*}
 The \emph {$(k+1)$-st $($geometric$)$  prolongation} of the bundle $P^0$ is the bundle $P^{k+1}$ over $P^k$ def\/ined by
 \begin{gather*}
 P^{k+1}=\big\{(\lambda_k, \mathcal H_{k+1}):\lambda_k\in P^k, C_{\mathcal H_{k+1}}^k\in \mathcal N_k\big\}.
 \end{gather*}
Equivalently,
\begin{gather*}
 P^{k+1}=\big\{(\lambda,\vf^{\mathcal H_{k+1}}):\lambda_k\in P^k, C_{\mathcal H_{k+1}}^k\in \mathcal N_k\big\}.
 \end{gather*}
It is a principal bundle with the abelian structure group $G^{k+1}$ of all maps $A\in \displaystyle{\bigoplus_{i\leq k}}{\rm Hom}(\mg^i$, $\mg^{i}\oplus\mg^{i+k+1}) $ such that
\begin{gather*}
A|_{\mg^i}=\begin{cases}{\rm Id}_{\mg^i}+T_i  &\text{ if } i<0,\\
{\rm Id}_{\mg^i} &\text{ if } 0\leq i\leq k,
\end{cases}
\end{gather*}
where $T_i\in {\rm Hom}(\mg^i,\mg^{i+k+1})$ and $(T_{-\mu},\ldots,T_{-1})\in \mg^{k+1}$. The right action $R_A^{k+1}$ of $A\in G^{k+1}$ on a f\/iber of $P^{k+1}$  is def\/ined by $R_A^{k+1} (\vf^{\mathcal H_{k+1}}) =\vf^{\mathcal H_{k+1}}\circ A$. Obviously, $G^{k+1}$ is an abelian group of dimension equal to $\dim \mg^{k+1}$. It is easy to see that the bundle $P^{k+1}$
is constructed so that the Properties $1$--$4$, formulated in the beginning of the present section, hold for $l=k+1$ as well.

Finally, assume that there exists $\bar l\geq 0$ such that $\mg^{\bar l}\neq 0$ but $\mg^{\bar l+1}= 0$. Since the symbol $\mathfrak m$ is fundamental, it follows that $\mg ^l=0$ for all $l>\bar l$. Hence, for all $l>\bar l$ the f\/iber of $P^l$ over a~point $\lambda_{l-1}\in P^{l-1}$ is a single point belonging to
$\displaystyle{\bigoplus_{i=-\mu}^{l-1} \text{Hom}\bigl(\mg^i,
D_{l-1}^i(\lambda_{l-1})/ D_{l-1}^{i+l+1}(\lambda_{l-1})\bigr)},$
where, as before, $\mu$ is the degree of nonholonomy of the distribution $D$.
Moreover, by our assumption, $D_l^i=0$ if $l\geq\bar l$ and $i\geq \bar l$. Therefore, if $l=\bar l+\mu$, then $i+l+1>\bar l$ for $i\geq -\mu$ and the f\/iber of $P^l$ over $P^l$ is an element of ${\rm Hom} \Big(\displaystyle{\bigoplus_{i=-\mu}^{l-1}}\mg^i, T_{\lambda_{l-1}}P^{l-1}\Big)$. In other words, $P^{\bar l+\mu}$ def\/ines a~canonical frame on $P^{\bar l+\mu-1}$.
But all bundles $P^l$ with $l\geq \bar l$ are identif\/ied one with each other by the canonical projections (which are dif\/feomorphisms in that case). As a conclusion we get an alternative proof of the main result of the Tanaka paper \cite{tan}:

 \begin{theorem}
 \label{main}
If the $(\bar l+1)$-st algebraic prolongation of the graded Lie algebra $\mathfrak m\oplus\mg^0$ is equal to zero then for any structure $P^0$ of constant
type $(\mathfrak m,\mg^0)$ there exists a canonical frame on the $\bar l$-th geometric prolongation $P^{\bar l}$ of $P^0$.
\end{theorem}

The power of Theorem \ref{main} is that it reduces the question of existence of a canonical frame for a~structure of constant type
$(\mathfrak m,\mg^0)$ to the calculation of the universal algebraic prolongation of the algebra $\mathfrak m\oplus \mg^0$. But the latter is pure Linear Algebra: each consecutive algebraic prolongation is determined by solving the system of linear equations given by~\eqref{mgk}. Let us demonstrate this algebraic prolongation procedure in the case of the equivalence of second order ordinary dif\/ferential equations with respect to the group of point transformations (see Example~5 in Subsection~\ref{section1.3}).
The result of this prolongation is very well known using the structure theory of simple Lie algebras (see discussions below), but this is one of the few nontrivial examples, where explicit calculations of algebraic prolongation can be written down in detail within one and a half pages.

{\bf Continuation of Example 5.} Recall that our geometric structure here is a contact distribution $D$ on a $3$-dimensional manifold endowed with two distinguished transversal line sub-distributions. The symbol
of $D$ is isomorphic to the $3$-dimensional Heisenberg algebra $\eta_3$ with grading $\mg^{-1}\oplus\mg^{-2}$, where $\mg^{-2}$ is the center of
$\eta_3$. Besides, the plane $\mg^{-1}$ is endowed with two distinguished transversal lines $\ell_1$ and $\ell_2$.  Let $X_1$ and $X_2$ be vectors spanning $\ell_1$ and $\ell_2$, respectively, and let $X_3=[X_1,X_2]$. Let $g^0$ be the algebra of all derivations on $\eta_3$ preserving the grading and the lines $\ell_1$ and $\ell_2$. Then
\begin{gather*}
\mg^0=\text{span}\big\{\Lambda_1^0,\Lambda_2^0\big\},
\end{gather*}
where
\begin{gather}
\label{2ordlam}
\Lambda_1^0(X_1)=X_1,\qquad  \Lambda_1^0(X_2)=X_2,\qquad
\Lambda_2^0(X_1)=X_1,\qquad \Lambda_2^0(X_2)=-X_2.
\end{gather}
Using the fact that $\Lambda_i^0$ is a derivation, we also have
\begin{gather}
\label{2ordlamx3}
\Lambda_1^0(X_3)=2X_3,\qquad \Lambda_2^0(X_3)=0.
\end{gather}

 {\bf a) Calculation of $\boldsymbol{\mg^{1}}$.}
 Given $\delta^1\in {\rm Hom}(\mg^{-1},\mg^0)\oplus {\rm Hom}(\mg^{-2},\mg^{-1})$ we have
 \begin{gather}
 \label{2ordg1}
 \delta^1(X_1)=\alpha_{11}\Lambda_1^0+\alpha_{12}\Lambda_2^0,
\qquad \delta^1(X_2)=\alpha_{21}\Lambda_1^0+\alpha_{22}\Lambda_2^0
 \end{gather}
 for some $\alpha_{ij}$, $1\leq i,j\leq 2$. If $\delta^1\in \mg^1$, then \eqref{2ordlam} yields
 \begin{gather}
 \label{2orddelx3}
\delta^1(X_3)=[\delta^1(X_1), X_2]+[X_1,\delta^1(X_2)]\\
\phantom{\delta^1(X_3)}{} =\bigl(\alpha_{11}\Lambda_1^0\!+\alpha_{12}\Lambda_2^0\bigr)(X_2)\!-\!
 \bigl(\alpha_{21}\Lambda_1^0\!+\alpha_{22}\Lambda_2^0\bigr)(X_1)=
 -(\alpha_{21}\!+\alpha_{22})X_1\!+(\alpha_{11}\!-\alpha_{12})X_2.\!
 \nonumber
 \end{gather}
 Using \eqref{2ordlamx3} and \eqref{2orddelx3}, we have
 \begin{gather*}
  0=\delta^1([X_1,X_3])=[\delta^1(X_1), X_3]+[X_1,\delta^1(X_3)]\\
  \phantom{0}{} = (\alpha_{11}\Lambda_1^0+\alpha_{12}\Lambda_2^0\bigr)(X_3)+(\alpha_{11}-\alpha_{12})X_3=
 (3\alpha_{11}-\alpha_{12})X_3,
 \end{gather*}
 which implies that
 $\alpha_{12}=
 3\alpha_{11}$.
  In the same way, from the identities $0=\delta^1[X_2,X_3]=[\delta^1(X_2), X_3]+[X_2,\delta^1(X_3)]$ one obtains easily that
  $\alpha_{22}=-3\alpha_{21}$.
   This completes the verif\/ication of conditions for $\delta^1$ to be in $\mg^{1}$. Hence
 \begin{gather*}
\mg^1=\text{span}\big\{\Lambda_1^1,\Lambda_2^1\big\},
\end{gather*}
where $\Lambda_1^1,\Lambda_2^1 \in {\rm Hom}(\mg^{-1},\mg^0)\oplus {\rm Hom}(\mg^{-2},\mg^{-1})$ such that
\begin{alignat}{4}
 &\Lambda_1^1(X_1)=\Lambda_1^0+3\Lambda_2^0,\qquad & & \Lambda_1^1(X_2)=0,\qquad & &\Lambda_1^1(X_3)=-2X_2,& \nonumber\\
 &\Lambda_2^1(X_1)=0,\qquad & & \Lambda_2^1(X_2)=\Lambda_1^0-3\Lambda_2^0,\qquad  && \Lambda_2^1(X_3)=2X_1.& \label{2ordlam1}
\end{alignat}
($\Lambda_1^1$ corresponds to $\delta^1$ as in \eqref{2ordg1} and \eqref{2orddelx3} with $\alpha_{11}=1$ and $\alpha_{21}=0$,
while $\Lambda_2^1$ corresponds to $\delta^1$ with $\alpha_{11}=0$ and $\alpha_{21}=1$ in the same formulas).

{\bf b) Calculation of $\boldsymbol{\mg^{2}}$.}
 Given $\delta^2\in {\rm Hom}(\mg^{-1},\mg^1)\oplus {\rm Hom}(\mg^{-2},\mg^{0})$ we have
 \begin{gather*}
 \delta^2(X_1)=\beta_{11}\Lambda_1^1+\beta_{12}\Lambda_2^1,
\qquad \delta^2(X_2)=\beta_{21}\Lambda_1^1+\beta_{22}\Lambda_2^1
 \end{gather*}
 for some $\beta_{ij}$, $1\leq i,j\leq 2$. If $\delta^2\in \mg^2$, then \eqref{2ordlam1} implies
\begin{gather}\label{2orddelx32}
 \delta^2(X_3)=[\delta^2(X_1), X_2]+[X_1,\delta^2(X_2)] \\
 \phantom{\delta^2(X_3)}{} =\bigl(\beta_{11}\Lambda_1^1+\beta_{12}\Lambda_2^1\bigr)(X_2)-
 \bigl(\beta_{21}\Lambda_1^1+\beta_{22}\Lambda_2^1\bigr)(X_1)=
 (\beta_{12}\!-\beta_{21})\Lambda_1^0\!-3(\beta_{12}\!+\beta_{21})\Lambda_2^0.\nonumber
 \end{gather}
 Using \eqref{2ordlam1} and \eqref{2orddelx32}, we have
 \begin{gather*}
  0=\delta^2([X_1,X_3])=[\delta^2(X_1), X_3]+[X_1,\delta^2(X_3)] \\
\phantom{0}{}= (\beta_{11}\Lambda_1^1+\beta_{12}\Lambda_2^1\bigr)(X_3)-
(\beta_{12}-\beta_{21})\Lambda_1^0(X_1)+3(\beta_{12}+\beta_{21})\Lambda_2^0(X_1)\\
\phantom{0}{}  =4(\beta_{12}+\beta_{21})X_1-2\beta_{11}X_2,
\end{gather*}
 which implies that
 \begin{gather}
 \label{2ordcond11}
 \beta_{11}=0, \qquad \beta_{21}=-\beta_{12}.
 \end{gather}
  Similarly, the identities $0=\delta^2[X_2,X_3]=[\delta^2(X_2), X_3]+[X_2,\delta^2(X_3)]$ implies
  $\beta_{22}=0$ in addition to \eqref{2ordcond11}.
  This completes the verif\/ications of conditions for $\delta^2$ to be in $\mg^{2}$. Hence
 \begin{gather*}
\mg^2=\text{span}\{\Lambda\},
\end{gather*}
where $\Lambda \in {\rm Hom}(\mg^{-1},\mg^1)\oplus {\rm Hom}(\mg^{-2},\mg^{0})$ is def\/ined by
\begin{gather}
\label{2ordlam2}
\Lambda(X_1)=\Lambda_2^1, \qquad \Lambda(X_2)=-\Lambda_1^1,\qquad \Lambda(X_3)=2\Lambda_1^0.
\end{gather}

{\bf c) Calculation of $\boldsymbol{\mg^{3}}$.}
Given $\delta^3\in {\rm Hom}(\mg^{-1},\mg^2)\oplus {\rm Hom}(\mg^{-2},\mg^{1})$ such that
 \begin{gather*}
 \delta^3(X_1)=\gamma_1\Lambda,
\qquad \delta^2(X_2)=\gamma_2\Lambda.
 \end{gather*}
 for some $\gamma_{i}$, $i=1,2$. If $\delta^3\in \mg^2$, then \eqref{2ordlam2} implies
 $\delta^3(X_1)=-\gamma_1\Lambda_1^1-\gamma_2\Lambda_2^1$. From the identities $0=\delta^3[X_1,X_3]=[\delta^3(X_1), X_3]+[X_1,\delta^3(X_3)]$ it follows easily that $\gamma_1=0$. From the identities $0=\delta^3[X_2,X_3]=[\delta^3(X_2), X_3]+[X_2,\delta^3(X_3)]$ it follows easily that  $\gamma_2=0$. Hence,
\[
\mg^3=0.
\]

  Thus the algebraic universal prolongation $\mg(\eta_3,\mg^0)=\mg^{-2}\oplus\mg^{-1}\oplus\mg^{0}\oplus\mg^{1}\oplus\mg^{2}$ is  $8$-dimensional. Therefore, f\/ixing the normalization conditions at each step, one can construct the f\/irst and the second
 geometric prolongations $P^1$ and $P^2$, and \emph{for any contact distribution $D$ on $3$-dimensional manifold endowed with two distinguished transversal line sub-distributions there is a canonical frame on the $8$-dimensional bundle $P^2$}.

 Let us look at the algebra $\mg(\eta_3,\mg^0)$ in more detail.
 Applying \eqref{posbr} inductively, we see that all nonzero brackets of elements $\Lambda_1^0$, $\Lambda_2^0$, $\Lambda_1^1$, $\Lambda_2^1$, $\Lambda$ (spanning the subalgebra of elements with nonzero weights of $\mg(\eta_3,\mg^0))$ are as follows:
 \[
 [\Lambda_1^1,\Lambda_1^0]=\Lambda_1^1,\!\!\qquad  [\Lambda_1^1,\Lambda_2^0]=\Lambda_1^1,\!\!\qquad
 [\Lambda_2^1,\Lambda_1^0]=\Lambda_1^1,\!\!\qquad [\Lambda_2^1,\Lambda_1^0]=-\Lambda_2^1,\!\!\qquad
[\Lambda_1^1,\Lambda_2^1]=2\Lambda.\!
\]
Considering all products in $\mg(\eta_3,\mg^0)$ it is not hard to see that $\mg(\eta_3,\mg^0)$ is isomorphic to $\mathfrak{sl}(3,\mathbb R)$. Indeed, if we denote by $E_{ij}$ the $3\times 3$-matrix such that its $(i,j)$ entry is equal to $1$ and all other entries vanish, then the following mapping  is an isomorphism of algebras $\mg(\eta_3,\mg^0)$ and $\mathfrak{sl}(3,\mathbb R)$:
\begin{gather*}
 X_1\mapsto E_{12},\qquad  X_2\mapsto E_{23},\qquad  X_3\mapsto E_{13},\\
 \Lambda_1^1\mapsto -2E_{21},\qquad \Lambda_2^1\mapsto-2E_{23},\qquad \Lambda\mapsto -2E_{23},\\
\Lambda_1^0+3\Lambda_2^0\mapsto 2(E_{11}-E_{22}),\qquad\Lambda_1^0-3\Lambda_2^0\mapsto 2(E_{22}-E_{33}).
\end{gather*}
As a matter of fact, here we are in the situation, when $\mathfrak m\oplus\mg^0$ is a subalgebra of elements of nonnegative weights
(a parabolic subalgbebra) of a graded simple Lie algebra (in the considered case $\eta_3\oplus\mg^0$ is a Borel subalgebra of $\mathfrak{sl}(3, \mathbb R)$). It was shown in \cite{yam} that, except for a few cases, the algebraic universal prolongation of a parabolic subalgebra of a simple Lie algebra is isomorphic to this simple Lie algebra.

This result can be applied also to the algebra $\mathfrak m_{(2,5)}\oplus\mg^0(\mathfrak m_{(2,5)})$, corresponding to maximally nonholonomic rank $2$ distributions in $\mathbb R^5$ (see Example 2 above).
In this case $\mathfrak m_{(2,5)}\oplus\mg^0(\mathfrak m_{(2,5)})$ is the subalgebra of elements of non-negative degree in the exceptional Lie algebra $G_2$, graded according to the coef\/f\/icient of the short simple root.
So, according to~\cite{yam}, the algebraic universal prolongation of $\mathfrak m_{(2,5)}\oplus\mg^0(\mathfrak m_{(2,5)})$ is isomorphic to $G_2=\displaystyle{\bigoplus_{i=-3}^3}\mg^i$. This together with Theorem~\ref{main} implies that \emph{to any
maximally nonholonomic rank $2$ distribution in $\mathbb R^5$ one can assign a canonical frame on the bundle $P^3$
 of dimension equal to $\dim G_2=14$}. Note that this statement is still weaker than what Cartan proved in \cite{cartan}.
  Indeed, Cartan provides explicit expressions for the coframe and f\/inds the complete system of invariants, while Theorem~\ref{main} is only the existence statement.

 Finally note that the construction of the bundles $P^k$ (and therefore of the canonical frame) depends on the choice of the normalization conditions given by spaces $\mathcal N_k$, as in \eqref{normsplitk}.
 Under additional assumptions on the algebra $\mg(\mathfrak m,\mg^0)$ (for example, semisimplicity or existence of a~special bilinear form) the spaces  $\mathcal N_k$ themselves can be taken in a canonical way at each step of the prolongation procedure. This allows to construct canonical frames satisfying additional nice properties.

 In particular, in another fundamental paper of Tanaka~\cite{tan2}, it was shown that if the algebraic universal prolongation  $\mg(\mathfrak m,\mg^0)$ is a semisimple Lie algebra, then the so-called $\mg(\mathfrak m,\mg^0)$-valued \emph{normal} Cartan connection can be associated with a structure of type $(\mathfrak m,\mg^0)$. Roughly \mbox{speaking}, a Cartan connection gives the canonical frame which is compatible in a natural way with the whole algebra  $\mg(\mathfrak m,\mg^0)$.
This is a generalization of Cartan's results \cite{cartan} on maximally nonholonomic rank 2 distributions in $\mathbb R^5$.
 Further,  T.~Morimoto \cite{mori} gave a general criterion (in terms of the algebra $\mg(\mathfrak m,\mg^0)$) for the existence of the normal Cartan connection for structures of type~$(\mathfrak m,\mg^0)$.

 All these developments are far beyond of the goals of the present note, so we do not want to address them in more detail here, referring the reader to the original papers.

\subsection*{Acknowledgements}

I would like to thank my colleagues Joseph Landsberg, Colleen Robles, and Dennis The for encouraging me to give these lectures, and  Boris Doubrov for stimulating discussions. I am also very grateful to the anonymous referees for numerous useful suggestions and corrections.

\pdfbookmark[1]{References}{ref}
\LastPageEnding

\end{document}